\documentclass[12pt,a4paper]{article}

\usepackage[left=2cm,right=2cm,    top=2cm,bottom=2cm,bindingoffset=0cm]{geometry}

\usepackage{subfig}
\usepackage{flushend}
\usepackage{indentfirst}
\usepackage{graphics}
\usepackage{amsmath}
\usepackage{epstopdf}
\usepackage{float}
\usepackage{amssymb}
\usepackage[font=footnotesize,labelfont=bf]{caption}
\usepackage[labelsep=period]{caption}

\newcommand{\be}{\begin{equation}}
\newcommand{\ee}{\end{equation}}

\newtheorem{thm}{Theorem}[section]
\newtheorem{lem}{Lemma}[section]
\newtheorem{nas}{Corollary}[section]

\newtheorem{zau}{Remark}[section]

\newtheorem{ozn}{Definition}[section]
\newtheorem{exm}{Example}[section]

\begin{document}

\title{\textbf{Robust extrapolation problem for stochastic sequences with stationary increments }}

\date{}

\maketitle

\noindent Columbia International Publishing\\
Contemporary Mathematics and Statistics
(2013) Vol. 1 No. 3 pp. 123-150\\
doi:10.7726/cms.2013.1009

\vspace{20pt}

\author{\textbf{Maksym Luz}, \textbf{Mykhailo Moklyachuk}$^*$, \\\\
 {Department of Probability Theory, Statistics and Actuarial
Mathematics, \\
Taras Shevchenko National University of Kyiv, Kyiv 01601, Ukraine}\\
$^{*}$Corresponding Author: Moklyachuk@gmail.com}\\\\\\

\noindent \textbf{Abstract.} \hspace{2pt}
The problem of optimal estimation of functionals $A\xi =\sum\nolimits_{k=0}^{\infty }{}a(k)\xi (k)$ and ${{A}_{N}}\xi =\sum\nolimits_{k=0}^{N}{}a(k)\xi (k)$
which depend on the unknown values of stochastic sequence $\xi (k)$ with stationary $n$th increments is considered.
Estimates are based on observations of the sequence $\xi (m)$ at points of time $m=-1,-2,\ldots$.
Formulas for calculating the value of the mean square error and the spectral characteristic of the optimal linear estimates of the functionals
are derived in the case where spectral density of the sequence is exactly known. Formulas that determine the least favorable spectral densities
and minimax (robust) spectral characteristic of the optimal linear estimates of the functionals are proposed in the case
where the spectral density of the sequence is not known but a set of admissible spectral densities is given.\\

\noindent \textbf{Keywords:} \hspace{2pt} Stochastic sequence with stationary increments; Minimax-robust estimate; Mean square error; Least favorable spectral density; Minimax-robust spectral characteristic.\\

\noindent \textbf{Mathematics Subject Classification:} \hspace{2pt}
 Primary: 60G10, 60G25, 60G35, Secondary: 62M20, 93E10

\section{{Introduction}}

Stochastic processes with $n$th stationary increments ${{\xi }^{(n)}}(t,\mu )$, $t,\mu \in R$, were introduced by Yaglom (1955). He described the main properties of these processes, found the spectral representation of stationary increments and solved the extrapolation problem for processes with stationary increments. Further results for such stochastic processes were presented by Pinsker (1955), Yaglom and Pinsker (1954). See Yaglom (1987a, 1987b) for more relative results and references.

The mean square optimal estimation problems for stochastic processes with th stationary increments are natural generalization of the linear extrapolation, interpolation and filtering problems for stationary stochastic processes.

Traditional methods of solution of the linear extrapolation, interpolation and filtering problems for stationary stochastic processes were developed by A.N. Kolmogorov, N.Wiener, A.M.Yaglom (see, for example, selected works of Kolmogorov (1992), survey article by Kailath (1974), books by Rozanov (1967), Wiener (1966), Yaglom (1987a, 1987b)). These methods are based on the assumption that the spectral density of the process is known.

In practice, however, it is impossible to have complete information on the spectral density in most cases. To solve the problem one finds parametric or nonparametric estimates of the unknown spectral density or selects a density by other reasoning. Then the classical estimation method is applied provided that the estimated or selected density is the true one. This procedure can result in significant increasing of the value of error as Vastola and Poor (1983) have demonstrated with the help of some examples. This is a reason to search estimates which are optimal for all densities from a certain class of admissible spectral densities. These estimates are called minimax since they minimize the maximal value of the error. A survey of results in minimax (robust) methods of data processing can be found in the paper by Kassam and Poor (1985). The paper by Ulf Grenander (1957) should be marked as the first one where the minimax extrapolation problem for stationary processes was formulated and solved. Franke and Poor (1984), Franke (1985) investigated the minimax extrapolation and filtering problems for stationary sequences with the help of convex optimization methods. This approach makes it possible to find equations that determine the least favorable spectral densities for various classes of admissible densities. For more details see, for example, books by Moklyachuk (2008), Moklyachuk and Masyutka (2012). In papers by Moklyachuk (1994-2008) the minimax approach was applied to extrapolation, interpolation and filtering problems for functionals which depend on the unknown values of stationary processes and sequences. Methods of solution the minimax-robust estimation problems for vector-valued stationary sequences and processes were developed by Moklyachuk and Masyutka (2006-2011). The minimax-robust estimation problems (extrapolation, interpolation and filtering) for linear functionals which depend on unknown values of periodically correlated stochastic processes were investigated by Dubovets'ka and Moklyachuk (2012-2013). Luz and Moklyachuk (2012a, 2012b) investigated the minimax interpolation problem for the linear functional ${{A}_{N}}\xi =\sum\nolimits_{k=0}^{N}{}a(k)\xi (k)$ which depends unknown values of a stochastic sequence $\xi (m)$ with stationary $n$th increments from observations of the sequence at points $\mathbb{Z}\backslash \{0,1,\ldots ,N\}$.

In this article we focus on the mean square optimal estimates of the functionals
\begin{equation}\label{eq01}
A\xi =\sum\limits_{k=0}^{\infty }{}(k)\xi (k),\quad {{A}_{N}}\xi =\sum\limits_{k=0}^{N}{}a(k)\xi (k)
\end{equation}
 which depend on the unknown values of a stochastic sequence $\xi (k)$ with stationary $n$th increments. Estimates are based on observations of the sequence $\xi (m)$ at points $m=-1,-2,\ldots $. The estimation problem for sequences with stationary increments is solved in the case of spectral certainty where the spectral density of the sequence is exactly known as well as in the case of spectral uncertainty where the spectral density of the sequence is not known but a set of admissible spectral densities is given. Formulas are derived for computing the value of the mean-square error and the spectral characteristic of the optimal linear estimates of functionals $A\xi $ and ${{A}_{N}}\xi $ in the case of spectral certainty. Formulas that determine the least favorable spectral densities and the minimax (robust) spectral characteristic of the optimal linear estimates of the functionals are proposed in the case of spectral uncertainty for concrete classes of admissible spectral densities.

\section{{ Stationary stochastic increment sequence. Spectral representation}}

\begin{ozn}\label{def2.1}
For a given stochastic sequence $\{\xi (m),m\in \mathbb{Z}\}$ a sequence
\begin{equation}\label{eq02}
	{{\xi }^{(n)}}(m,\mu )=(1-{{B}_{\mu }}{{)}^{n}}\xi (m)=\sum\limits_{l=0}^{n}{}{{(-1)}^{l}}C_{n}^{l}\xi (m-l\mu ),
\end{equation}
where ${{B}_{\mu }}$ is a backward shift operator with step $\mu \in \mathbb{Z}$, such that ${{B}_{\mu }}\xi (m)=\xi (m-\mu )$, is called the stochastic $n$th increment sequence with step $\mu \in \mathbb{Z}$.
\end{ozn}

For the stochastic $n$th increment sequence ${{\xi }^{(n)}}(m,\mu )$ the following relations hold true:
\begin{equation}\label{eq03}
	{{\xi }^{(n)}}(m,-\mu )=(-{{1)}^{n}}{{\xi }^{(n)}}(m+n\mu ,\mu ),
\end{equation}
\begin{equation}\label{eq04}
	{{\xi }^{(n)}}(m,k\mu )=\sum\limits_{l=0}^{(k-1)n}{}{{A}_{l}}{{\xi }^{(n)}}(m-l\mu ,\mu ),\quad k\in N,
\end{equation}
where coefficients $\{{{A}_{l}},l=0,1,2,\ldots ,(k-1)n\}$ are determined by the representation
$${{(1+x+\ldots +{{x}^{k-1}})}^{n}}=\sum\limits_{l=0}^{(k-1)n}{}{{A}_{l}}{{x}^{l}}.$$

\begin{ozn}\label{def2.2}
 The stochastic $n$th increment sequence ${{\xi }^{(n)}}(m,\mu )$ generated by stochastic sequence $\{\xi (m),m\in \mathbb{Z}\}$ is wide sense stationary if the mathematical expectations
	$$\text{E}{{\xi }^{(n)}}({{m}_{0}},\mu )={{c}^{(n)}}(\mu )$$
and
	$$\text{E}{{\xi }^{(n)}}({{m}_{0}}+m,{{\mu }_{1}}){{\xi }^{(n)}}({{m}_{0}},{{\mu }_{2}})={{D}^{(n)}}(m,{{\mu }_{1}},{{\mu }_{2}})$$
exist for all ${{m}_{0}},\mu ,m,{{\mu }_{1}},{{\mu }_{2}}$ and do not depend on ${{m}_{0}}$.
The function ${{c}^{(n)}}(\mu )$ is called the mean value of the $n$th increment sequence and the function ${{D}^{(n)}}(m,{{\mu }_{1}},{{\mu }_{2}})$ is called the structural function of the stationary $n$th increment sequence (or the structural function of $n$th order of the stochastic sequence $\{\xi (m),m\in \mathbb{Z}\}$).
\\
The stochastic sequence $\{\xi (m),m\in \mathbb{Z}\}$ which determines the stationary $n$th increment sequence ${{\xi }^{(n)}}(m,\mu )$ by formula (2) is called sequence with stationary $n$th increments.
\end{ozn}

\begin{thm}\label{thm2.1}
 The mean value ${{c}^{(n)}}(\mu )$ and the structural function ${{D}^{(n)}}(m,{{\mu }_{1}},{{\mu }_{2}})$ of the stochastic stationary $n$th increment sequence ${{\xi }^{(n)}}(m,\mu )$ can be represented in the following forms
\begin{equation}\label{eq05}
{{c}^{(n)}}(\mu )=c{{\mu }^{n}},
\end{equation}
\begin{equation}\label{eq06}
{{D}^{(n)}}(m,{{\mu }_{1}},{{\mu }_{2}})=\int\limits_{-\pi }^{\pi }{{{e}^{i\lambda m}}}{{(1-{{e}^{-i{{\mu }_{1}}\lambda }})}^{n}}{{(1-{{e}^{i{{\mu }_{2}}\lambda }})}^{n}}\frac{1}{{{\lambda }^{2n}}}dF(\lambda ),
\end{equation}
where $c$ is a constant, $F(\lambda )$ is a left-continuous nondecreasing bounded function with $F(-\pi )=0.$ The constant $c$ and the function $F(\lambda )$ are determined uniquely by the increment sequence ${{\xi }^{(n)}}(m,\mu )$.
\\
From the other hand, a function ${{c}^{(n)}}(\mu )$ which has the form $(5)$ with a constant $c$ and a function ${{D}^{(n)}}(m,{{\mu }_{1}},{{\mu }_{2}})$ which has the form $(6)$ with a function $F(\lambda )$ which satisfies the indicated conditions are the mean value and the structural function of some stationary $n$th increment sequence ${{\xi }^{(n)}}(m,\mu )$.
\end{thm}

Using representation $(6)$ of the structural function of a stationary $n$th increment sequence ${{\xi }^{(n)}}(m,\mu )$ and the Karhunen theorem (see Karhunen (1947)), we get the following spectral representation of the stationary $n$th increment sequence ${{\xi }^{(n)}}(m,\mu )$:
\begin{equation}\label{eq07}
{{\xi }^{(n)}}(m,\mu )=\int\limits_{-\pi }^{\pi }{{{e}^{i\lambda m}}}{{(1-{{e}^{-i\mu \lambda }})}^{n}}\frac{1}{{{(i\lambda )}^{n}}}dZ(\lambda ),
\end{equation}
where $Z(\lambda )$ is an orthogonal stochastic measure ?n $[-\pi ,\pi )$ connected with the spectral function $F(\lambda )$ by the relation
\begin{equation}\label{eq08}
\text{E}Z({{A}_{1}})\overline{Z({{A}_{2}})}=F({{A}_{1}}\cap {{A}_{2}}).
\end{equation}
Denote by $H({{\xi }^{(n)}})$ the Hilbert space generated by all elements $\{{{\xi }^{(n)}}(m,\mu ):m,\mu \in \mathbb{Z}\}$ in the space $H={{L}_{2}}(\Omega ,F,P)$ and let ${{H}^{t}}({{\xi }^{(n)}})$, $t\in \mathbb{Z}$, be the subspace of $H({{\xi }^{(n)}})$ generated by elements $\{{{\xi }^{(n)}}(m,\mu ):m\le t,\mu >0\}$. Let
$$S({{\xi }^{(n)}})=\bigcap\limits_{t\in }{}{{H}^{t}}({{\xi }^{(n)}}).$$
Since the space $S({{\xi }^{(n)}})$ is a subspace in the Hilbert space $H({{\xi }^{(n)}})$, the space $H({{\xi }^{(n)}})$ admits the decomposition
$$H({{\xi }^{(n)}})=S({{\xi }^{(n)}})\oplus R({{\xi }^{(n)}}),$$
where $R({{\xi }^{(n)}})$ is an orthogonal complement of the subspace $S({{\xi }^{(n)}})$ in the space $H({{\xi }^{(n)}})$.

\begin{ozn}\label{def2.3}
 A stationary $n$th increment sequence ${{\xi }^{(n)}}(m,\mu )$ is called regular if $H({{\xi }^{(n)}})=R({{\xi }^{(n)}})$. It is called singular if $H({{\xi }^{(n)}})=S({{\xi }^{(n)}})$.
\end{ozn}

\begin{thm}\label{thm2.2}
A wide-sense stationary stochastic increment sequence admits a unique representation in the form
\begin{equation}\label{eq09}
{{\xi }^{(n)}}(m,\mu )=\xi _{r}^{(n)}(m,\mu )+\xi _{s}^{(n)}(m,\mu ),
\end{equation}
where $\{\xi _{r}^{(n)}(m,\mu ):m\in \mathbb{Z}\}$ is a regular increment sequence and $\{\xi _{s}^{(n)}(m,\mu ):m\in \mathbb{Z}\}$ is a singular increment sequence. Moreover, the increment sequences $\xi _{r}^{(n)}(m,\mu )$ and $\xi _{s}^{(n)}(k,\mu )$ are orthogonal for all $m,k\in \mathbb{Z}$.
\end{thm}
Components of representation $(9)$ are constructed in the following way
$$\xi _{s}^{(n)}(m,\mu )=\text{E}[{{\xi }^{(n)}}(m,\mu )|S({{\xi }^{(n)}})],\quad \xi _{r}^{(n)}(m,\mu )={{\xi }^{(n)}}(m,\mu )-\xi _{s}^{(n)}(m,\mu ).$$

Let $\{{{\varepsilon }_{m}}:m\in \mathbb{Z}\}$ be a sequence of uncorrelated random variables with $\text{E}{{\varepsilon }_{m}}=0$ and $\text{D}\varepsilon _{m}^{2}=1$. Define the Hilbert space ${{H}^{t}}(\varepsilon )$ generated by elements $\{{{\varepsilon }_{m}}:m\le t\}$.

\begin{ozn}\label{def2.4}
   A sequence of uncorrelated random variables $\{{{\varepsilon }_{m}}:m\in \mathbb{Z}\}$ is called innovation sequence for a regular stationary $n$th increment sequence ${{\xi }^{(n)}}(m,\mu )$ if the condition ${{H}^{t}}({{\xi }^{(n)}})={{H}^{t}}(\varepsilon )$ holds true for all $t\in \mathbb{Z}$.
\end{ozn}

\begin{thm}\label{thm2.3}
  A stochastic stationary increment sequence ${{\xi }^{(n)}}(m,\mu )$ is regular if and only if there exists an innovation sequence $\{{{\varepsilon }_{m}}:m\in \mathbb{Z}\}$ and a sequence of complex functions $\{{{\varphi }^{(n)}}(k,\mu ):m\ge 0\}$ ,$\sum\nolimits_{k=0}^{\infty }{}|{{\varphi }^{(n)}}(k,\mu {{)|}^{2}}<\infty $, such that
\begin{equation}\label{eq10}
	{{\xi }^{(n)}}(m,\mu )=\sum\limits_{k=0}^{\infty }{}{{\varphi }^{(n)}}(k,\mu )\varepsilon (m-k).
\end{equation}
\end{thm}

Representation (10) is called canonical moving average representation of the stochastic stationary increment sequence ${{\xi }^{(n)}}(m,\mu )$.

\begin{nas} \label{corr2.1}
Corollary 2.1 A wide-sense stationary stochastic increment sequence admits a unique representation in the form
\begin{equation}\label{eq11}
	{{\xi }^{(n)}}(m,\mu )=\xi _{s}^{(n)}(m,\mu )+\sum\limits_{k=0}^{\infty }{}{{\varphi }^{(n)}}(k,\mu )\varepsilon (m-k),
\end{equation}
 where $\sum\nolimits_{k=0}^{\infty }{}|{{\varphi }^{(n)}}(k,\mu {{)|}^{2}}<\infty $ and $\{{{\varepsilon }_{m}}:m\in \mathbb{Z}\}$ is the innovation sequence.
\end{nas}

Let the stationary $n$th increment sequence ${{\xi }^{(n)}}(m,\mu )$ admit the canonical representation (10). In this case the spectral function $F(\lambda )$ of the stationary increment sequence ${{\xi }^{(n)}}(m,\mu )$ has a spectral density $f(\lambda )$ which admits the canonical factorization
\begin{equation}\label{eq12}
f(\lambda )=|\Phi ({{e}^{-i\lambda }}{{)|}^{2}},\quad \Phi (z)=\sum\limits_{k=0}^{\infty }{}\varphi (k){{z}^{k}},
\end{equation}
where the function $\Phi (z)=\sum\nolimits_{k=0}^{\infty }{}\varphi (k){{z}^{k}}$ has the convergence radius $r>1$ and does not have zeros in the unit disk $\{z:|z|\le 1\}$.
\noindent Let us define
$${{\Phi }_{\mu }}(z)=\sum\limits_{k=0}^{\infty }{}{{\varphi }^{(n)}}(k,\mu ){{z}^{k}}=\sum\limits_{k=0}^{\infty }{}{{\varphi }_{\mu }}(k){{z}^{k}},$$
where ${{\varphi }_{\mu }}(k)={{\varphi }^{(n)}}(k,\mu )$ are coefficients which determine the canonical representation $(10)$.
\noindent Then the following relation holds true
\begin{equation}\label{eq13}
	{{\left| {{\Phi }_{\mu }}({{e}^{-i\lambda }}) \right|}^{2}}=\frac{|1-{{e}^{-i\lambda \mu }}{{|}^{2n}}}{{{\lambda }^{2n}}}f(\lambda ).
\end{equation}

The one-sided moving average representation (10) and relation (13) are used for finding the optimal mean square estimate of the unknown values of a sequence with $n$th stationary increment.

\section{{Hilbert space projection method of extrapolation of linear functionals}}

Let $\{\xi (m),m\in \mathbb{Z}\}$ be a stochastic sequence which determines a stationary $n$th increment sequence ${{\xi }^{(n)}}(m,\mu )$ with an absolutely continuous spectral function $F(\lambda )$ which has spectral density $f(\lambda )$. Without loss of generality we will assume that the mean value of the increment sequence ${{\xi }^{(n)}}(m,\mu )$ is 0. Let the stationary increment sequence ${{\xi }^{(n)}}(m,\mu )$ admit the one-sided moving average representation $(10)$ and the spectral density $f(\lambda )$ admits the canonical factorization $(12)$.
Consider the case where the step $\mu >0$.

Suppose that observations of the sequence $\xi (m)$ at points $m=-1,-2,...$ are known. The problem is to find the mean square optimal linear estimates of functionals
$${{A}_{N}}\xi =\sum\nolimits_{k=0}^{N}{}a(k)\xi (k)$$, $$A\xi =\sum\nolimits_{k=0}^{\infty }{}a(k)\xi (k)$$ which depend on unknown values $\xi (m)$, $m\ge 0$ of the sequence $\xi (m)$.
From $(2)$ we can obtain the formal equation
\begin{equation}\label{eq14}
\xi (k)=\frac{1}{{{(1-{{B}_{\mu }})}^{n}}}{{\xi }^{(n)}}(k,\mu )=\sum\limits_{j=-\infty }^{k}{}{{d}_{\mu }}(k-j){{\xi }^{(n)}}(j,\mu ),
\end{equation}
where coefficients  $\{{{d}_{\mu }}(k):k\ge 0\}$ are determined by the relation
$$\sum\limits_{k=0}^{\infty }{}{{d}_{\mu }}(k){{x}^{k}}={{\left( \sum\limits_{j=0}^{\infty }{}{{x}^{\mu j}} \right)}^{n}}.$$
From $(2)$ and $(14)$ one can obtain the following relations
$$\sum\limits_{k=0}^{\infty }{}a(k)\xi (k)=-\sum\limits_{i=-\mu n}^{-1}{}v(i)\xi (i)+\sum\limits_{i=0}^{\infty }{}\left( \sum\limits_{k=i}^{\infty }{}a(k){{d}_{\mu }}(k-i) \right){{\xi }^{(n)}}(i,\mu ),$$
	$$\sum\limits_{k=0}^{\infty }{}{{b}_{\mu }}(k){{\xi }^{(n)}}(k,\mu )=\sum\limits_{i=-\mu n}^{-1}{}\xi (i)\sum\limits_{l={{\left[ -\frac{i}{\mu } \right]}^{'}}}^{n}{}{{(-1)}^{l}}C_{n}^{l}{{b}_{\mu }}(l\mu +i)+\sum\limits_{i=0}^{\infty }{}\xi (i)\sum\limits_{l=0}^{n}{}{{(-1)}^{l}}C_{n}^{l}{{b}_{\mu }}(l\mu +i),$$
where $[x{]}'$ denotes the least integer number among numbers which are greater or equal to $x$. Using these relations we obtain representation of the functional $A\xi $ as difference $A\xi =B\xi -V\xi $ of functionals, where
$$B\xi =\sum\limits_{k=0}^{\infty }{}{{b}_{\mu }}(k){{\xi }^{(n)}}(k,\mu ),\quad V\xi =\sum\limits_{k=-\mu n}^{-1}{}{{v}_{\mu }}(k)\xi (k),$$
\begin{equation}\label{eq15}
{{v}_{\mu }}(k)=\sum\limits_{l={{\left[ -\frac{k}{\mu } \right]}^{'}}}^{n}{}{{(-1)}^{l}}C_{n}^{l}{{b}_{\mu }}(l\mu +k),\quad k=-1,-2,\ldots ,-\mu n,
\end{equation}
\begin{equation}\label{eq16}
{{b}_{\mu }}(k)=\sum\limits_{m=k}^{\infty }{}a(m){{d}_{\mu }}(m-k)=({{D}^{\mu }}a{{)}_{k}},k\ge 0.
\end{equation}
Here ${{D}^{\mu }}$ is a linear operator in the space ${{\ell }_{2}}$ determined by elements $D_{k,j}^{\mu }={{d}_{\mu }}(j-k)$ if $0\le k\le j\ $, and $D_{k,j}^{\mu }=0$ if $j<k$; the vector $a=(a(0),a(1),a(2),\ldots {{)}^{T}}$.

We will suppose that the following restirictions on the coefficients $\{{{b}_{\mu }}(k):k\ge 0\}$ hold true
\begin{equation}\label{eq17}
\sum\limits_{k=0}^{\infty }{}|{{b}_{\mu }}(k)|<\infty ,\quad \sum\limits_{k=0}^{\infty }{}(k+1)|{{b}_{\mu }}(k{{)|}^{2}}<\infty.
\end{equation}
Under these conditions the functional $B\xi $ has the second moment and the operator ${{B}^{\mu }}$ defined below is compact. Since coefficients $a(k)$ and ${{b}_{\mu }}(k)$ are related by $(16)$, the following conditions hold true
\begin{equation}\label{eq18}
\sum\limits_{k=0}^{\infty }{}|({{D}^{\mu }}a{{)}_{k}}|<\infty ,\quad \sum\limits_{k=0}^{\infty }{}(k+1)|({{D}^{\mu }}a{{)}_{k}}{{|}^{2}}<\infty.
\end{equation}

Let $\hat{A}\xi $ denote the mean square optimal linear estimate of the functional $A\xi $ from observations of the sequence $\xi (m)$ at points $m=-1,-2,...$  and let $\hat{B}\xi $ denote the mean square optimal linear estimate of the functional $B\xi $ from observations of the stochastic $n$th increment sequence ${{\xi }^{(n)}}(m,\mu )$ at points $m=-1,-2,...$  . Let $\Delta (f,\hat{A}\xi ):=\text{E}|A\xi -\hat{A}\xi {{|}^{2}}$ denote the mean square error of the estimate $\hat{A}\xi $ and let $\Delta (f,\hat{B}\xi ):=\text{E}|B\xi -\hat{B}\xi {{|}^{2}}$ denote the mean square error of the estimate $\hat{B}\xi $. Since values of the sequence $\xi (m)$ are known for $m=-1,-2,\ldots ,-\mu n$, the following equality holds true
\begin{equation}\label{eq19}
\hat{A}\xi =\hat{B}\xi -V\xi.
\end{equation}
 From this relation we get

	$$\Delta (f,\hat{A}\xi )=\text{E}|A\xi -\hat{A}\xi {{|}^{2}}=\text{E}|A\xi +V\xi -\hat{B}\xi {{|}^{2}}=\text{E}|B\xi -\hat{B}\xi {{|}^{2}}=\Delta (f,\hat{B}\xi ).$$

Denote by $L_{2}^{0-}(f)$ the subspace of the Hilbert space ${{L}_{2}}(f)$ generated by the set of functions $\{{{e}^{i\lambda k}}{{(1-{{e}^{-i\lambda \mu }})}^{n}}\frac{1}{{{(i\lambda )}^{n}}}:k\le -1\}$. Every linear estimate $\hat{B}\xi $ of the functional $B\xi $ can be represented in the form
\begin{equation}\label{eq20}
\hat{B}\xi =\int\limits_{-\pi }^{\pi }{{{h}_{\mu }}(\lambda )(1-{{e}^{-i\lambda \mu }}{{)}^{n}}\frac{1}{{{(i\lambda )}^{n}}}}dZ(\lambda ),
\end{equation}
where ${{h}_{\mu }}(\lambda )$ is the spectral characteristic of the estimate $\hat{B}\xi $. The spectral characteristic of the optimal estimate provides the minimum value of the mean square error $\Delta (f,\hat{B}\xi )$.

With the help of the Hilbert space projection method proposed by Kolmogorov we can find formulas for calculation the mean square error and the spectral characteristic of the optimal linear estimate $\hat{B}\xi $ of the functional $B\xi $. Following the method we find that the the spectral characteristic ${{h}_{\mu }}(\lambda )$ of the optimal linear estimate is determined by the following conditions:

1) ${{h}_{\mu }}(\lambda )(1-{{e}^{-i\lambda \mu }}{{)}^{n}}\frac{1}{{{(i\lambda )}^{n}}}\in L_{2}^{0-}\left( f \right)$;

2) $({{B}^{\mu }}({{e}^{i\lambda }})-{{h}_{\mu }}(\lambda ))(1-{{e}^{-i\lambda \mu }}{{)}^{n}}\frac{1}{{{(i\lambda )}^{n}}}\bot L_{2}^{0-}\left( f \right)$, where

$${{B}^{\mu }}({{e}^{i\lambda }})=\sum\limits_{k=0}^{\infty }{}{{b}_{\mu }}(k){{e}^{i\lambda k}}.$$
From the second condition we obtain the following relation for every $k\le -1$

	$$\int\limits_{-\pi }^{\pi }{({{B}^{\mu }}({{e}^{i\lambda }})-{{h}_{\mu }}(\lambda ))|1-{{e}^{i\lambda \mu }}{{|}^{2n}}\frac{1}{{{\lambda }^{2n}}}{{e}^{-i\lambda k}}f(\lambda )d\lambda }=0.$$
 These relations are satisfied by the function
\begin{equation}\label{eq21}
{{h}_{\mu }}(\lambda )={{B}^{\mu }}({{e}^{i\lambda }})-{{r}_{\mu }}({{e}^{i\lambda }})\Phi _{\mu }^{-1}({{e}^{-i\lambda }}),
\end{equation}
$${{r}_{\mu }}({{e}^{i\lambda }})=\sum\limits_{j=0}^{\infty }{}\sum\limits_{m=0}^{\infty }{}{{b}_{\mu }}(m+j){{\varphi }_{\mu }}(m){{e}^{i\lambda j}}=\sum\limits_{j=0}^{\infty }{}{{({{B}^{\mu }}{{\varphi }_{\mu }})}_{j}}{{e}^{i\lambda j}},$$
 where ${{B}^{\mu }}$ is a linear symmetric operator in the space ${{\ell }_{2}}$ defined by the matrix with elements $B_{k,j}^{\mu }={{b}_{\mu }}(k+j)$, $k,j\ge 0$. ${{\varphi }_{\mu }}=({{\varphi }_{\mu }}(0),{{\varphi }_{\mu }}(1),{{\varphi }_{\mu }}(2),\ldots )$; ${{\varphi }_{\mu }}(k)$, $k\ge 0$, are coefficients which determine the moving average representation (10).

Note that under conditions $(17)$ the operator ${{B}^{\mu }}$ is compact.

To check condition 1) it is sufficient to show that the function ${{h}_{\mu }}(\lambda )\in L_{2}^{0-}$, where $L_{2}^{0-}$ is the closed linear subspace of the space ${{L}_{2}}(-\pi ,\pi )$ generated by the set of functions $\{{{e}^{i\lambda k}}:k\le -1\}$. Since $\Phi _{\mu }^{-1}({{e}^{-i\lambda }})\in L_{2}^{0-}$, we have
$${{h}_{\mu }}(\lambda )=({{B}_{\mu }}({{e}^{i\lambda }}){{\Phi }_{\mu }}({{e}^{-i\lambda }})-{{r}_{\mu }}({{e}^{i\lambda }}))\Phi _{\mu }^{-1}({{e}^{-i\lambda }})=$$
	$$=\Phi _{\mu }^{-1}({{e}^{-i\lambda }})\sum\limits_{j=-\infty }^{-1}{}\sum\limits_{m=-j}^{\infty }{}{{b}_{\mu }}(m+j){{\varphi }_{\mu }}(m){{e}^{i\lambda j}}\in L_{2}^{0-}.$$
Therefore the spectral characteristic ${{h}_{\mu }}(\lambda )=:{{h}_{\mu }}(f)$ of the optimal estimate $\hat{B}\xi $ of the functional $B\xi $ can be calculated by formula (21).

 The value of the mean square error $\Delta (f,\hat{B}\xi )$ can be calculated by the formula
\begin{equation}\label{eq22}
\Delta (f,\hat{B}\xi )=\frac{1}{2\pi }\int\limits_{-\pi }^{\pi }{|{{r}_{\mu }}({{e}^{i\lambda }}{{)|}^{2}}d\lambda }=||{{B}^{\mu }}{{\varphi }_{\mu }}{{||}^{2}}.
\end{equation}
Summarizing our reasoning we have the following theorem.

\begin{thm} \label{thm3.1}
  Let a stochastic sequence $\{\xi (m),m\in \mathbb{Z}\}$ determine a stationary stochastic $n$th increment sequence ${{\xi }^{(n)}}(m,\mu )$ with absolutely continuous spectral function $F(\lambda )$ and spectral density $f(\lambda )$ which admits the canonical factorization $(12)$. The optimal linear estimate $\hat{B}\xi $ of the functional $B\xi $ which depends on the unobserved values ${{\xi }^{(n)}}(m,\mu )$, $m=0,1,2,\ldots $, $\mu >0$, from observations of the sequence $\xi (m)$ at points $m=-1,-2,\ldots $, can be calculated by formula (20). The spectral characteristic ${{h}_{\mu }}(\lambda )$ of the optimal linear estimate $\hat{B}\xi $ can be calculated by formula $(21)$. The value of the mean square error $\Delta (f,\hat{B}\xi )$ can be calculated by formula $(22)$.
\end{thm}

Using Theorem $3.1$ and representation $(9)$, we can obtain the optimal estimate of an unobserved value of the sequence ${{\xi }^{(n)}}(m,\mu )$, $m\ge 0$, from observations of the sequence $\xi (k)$ at points $k=-1,-2,\ldots $ The singular component $\xi _{s}^{(n)}(k,\mu )$ of the sequence has errorless estimate. We will use formula $(21)$ to obtain the spectral characteristic ${{h}_{m,\mu }}(\lambda )$ of the optimal estimate ${{\hat{\xi }}^{(n)}}(m,\mu )$ of the regular component $\xi _{r}^{(n)}(k,\mu )$ of the sequence. Consider the vector ${{b}_{\mu }}$ with 1 on position $m$, $m\ge 0$, and 0 on other positions. It follows from the derived formulas that the spectral characteristic of the estimate
\begin{equation}\label{eq23}
{{\hat{\xi }}^{(n)}}(m,\mu )=\xi _{s}^{(n)}(k,\mu )+\int\limits_{-\pi }^{\pi }{{{h}_{m,\mu }}(\lambda )(1-{{e}^{-i\lambda \mu }}{{)}^{n}}\frac{1}{{{(i\lambda )}^{n}}}dZ(\lambda )}
\end{equation}
 can be calculated by the formula
\begin{equation}\label{eq24}
{{h}_{m,\mu }}(\lambda )={{e}^{i\lambda m}}-\Phi _{\mu }^{-1}({{e}^{-i\lambda }})\sum\limits_{k=0}^{m}{}{{\varphi }_{\mu }}(k){{e}^{-i\lambda k}}.
\end{equation}
The value of the mean square error can be calculated by the formula
\begin{equation}\label{eq25}
\Delta (f,{{\hat{\xi }}^{(n)}}(m,\mu ))=\frac{1}{2\pi }\int\limits_{-\pi }^{\pi }{{{\left| \sum\limits_{k=0}^{m}{}{{\varphi }_{\mu }}(k){{e}^{-i\lambda k}} \right|}^{2}}d\lambda }=\sum\limits_{k=0}^{m}{}|{{\varphi }_{\mu }}(k{{)|}^{2}}.
\end{equation}

The following statement holds true.

\begin{nas} \label{corr3.1}
The optimal linear estimate ${{\hat{\xi }}^{(n)}}(m,\mu )$ of the value of the increment sequence ${{\xi }^{(n)}}(m,\mu )$, $m\ge 0$, $\mu >0$, from observations of the sequence $\xi (k)$ at points $k=-1,-2,\ldots $ can be calculated by formula$(23)$. The spectral characteristic ${{h}_{m,\mu }}(\lambda )$ of the optimal linear estimate ${{\hat{\xi }}^{(n)}}(m,\mu )$ can be calculated by formula $(24)$. The value of mean square error $\Delta (f,{{\hat{\xi }}^{(n)}}(m,\mu ))$ of the optimal linear estimate can be calculated by formula $(25)$.
\end{nas}

Making use relation $(19)$ we can find the optimal estimate $\hat{A}\xi $ of the functional $A\xi $ from observations of the sequence $\xi (k)$ at points   $k=-1,-2,\ldots $. These estimate can be presented in the following form
\begin{equation}\label{eq26}
\hat{A}\xi =-\sum\limits_{k=-\mu n}^{-1}{}{{v}_{\mu }}(k)\xi (k)+\int\limits_{-\pi }^{\pi }{h_{\mu }^{(a)}(\lambda )(1-{{e}^{-i\lambda \mu }}{{)}^{n}}\frac{1}{{{(i\lambda )}^{n}}}dZ(\lambda )},
\end{equation}
where coefficients ${{v}_{\mu }}(k)$ for $k=-1,-2,\ldots ,-\mu n$ are defined by relation $(15)$. Using relationship $(16)$ between  coefficients $a(k)$ and ${{b}_{\mu }}(k)$, we obtain the following equation
$${{({{B}^{\mu }}{{\varphi }_{\mu }})}_{k}}=\sum\limits_{m=0}^{\infty }{}\sum\limits_{l=k}^{\infty }{}{{\varphi }_{\mu }}(m)a(m+l){{d}_{\mu }}(l-k)=({{D}^{\mu }}A{{\varphi }_{\mu }}{{)}_{k}},$$
where the linear operator $A$ is defined by coefficients $a(k)$, $k\ge 0$, in the following way: ${{(A)}_{k,j}}=a(k+j)$, $k,j\ge 0$. Thus the spectral characteristic and the value of the mean square error of the optimal estimate $\hat{A}\xi $ can be calculated by the formulas
\begin{equation}\label{eq27}
h_{\mu }^{(a)}(\lambda )=A({{e}^{i\lambda }})-r_{\mu }^{(a)}({{e}^{i\lambda }})\Phi _{\mu }^{-1}({{e}^{-i\lambda }}),
\end{equation}
\begin{equation}\label{eq28}
A({{e}^{i\lambda }})=\sum\limits_{k=0}^{\infty }{}{{({{D}^{\mu }}a)}_{k}}{{e}^{i\lambda k}},\quad r_{\mu }^{(a)}({{e}^{i\lambda }})=\sum\limits_{j=0}^{\infty }{}{{({{D}^{\mu }}A{{\varphi }_{\mu }})}_{j}}{{e}^{i\lambda j}}.
\end{equation}
\begin{equation}\label{eq29}
\Delta (f,\hat{A}\xi )=\frac{1}{2\pi }\int\limits_{-\pi }^{\pi }{|r_{\mu }^{(a)}({{e}^{i\lambda }}{{)|}^{2}}d\lambda }=||{{D}^{\mu }}A{{\varphi }_{\mu }}{{||}^{2}}.
\end{equation}

The following theorem holds true.

\begin{thm} \label{thm3.2}
 Let a stochastic sequence $\{\xi (m),m\in \mathbb{Z}\}$ determine a stationary stochastic $n$th increment sequence ${{\xi }^{(n)}}(m,\mu )$ with absolutely continuous spectral function $F(\lambda )$ and spectral density $f(\lambda )$ which admits the canonical factorization $(12)$. The optimal linear estimate $\hat{A}\xi $ of the functional $A\xi $ of unobserved values $\xi (m)$, $m=0,1,2,\ldots $, from observations of the sequence $\xi (m)$ at points $m=-1,-2,\ldots $, can be calculated by formula $(26)$. The spectral characteristic $h_{\mu }^{(a)}(\lambda )$ of the optimal linear estimate $\hat{A}\xi $ can be calculated by formula $(27)$. The value of the mean square error $\Delta (f,\hat{A}\xi )$ of the optimal linear estimate can be calculated by formula $(29)$.
\end{thm}

Consider now the problem of the mean square optimal estimation of the functional ${{A}_{N}}\xi $.

Using the derived formulas we can find the optimal estimate of the functional ${{A}_{N}}\xi $ in the form
\begin{equation}\label{eq30}
{{\hat{A}}_{N}}\xi =-\sum\limits_{k=-\mu n}^{-1}{}{{v}_{\mu ,N}}(k)\xi (k)+\int\limits_{-\pi }^{\pi }{h_{\mu ,N}^{(a)}(\lambda )(1-{{e}^{-i\lambda \mu }}{{)}^{n}}\frac{1}{{{(i\lambda )}^{n}}}dZ(\lambda )},
\end{equation}
where coefficients ${{v}_{\mu ,N}}(k)$, $k=-1,-2,\ldots ,-\mu n$, are calculated by formulas
$${{v}_{\mu ,N}}(k)=\sum\limits_{l=\left[ -\frac{k}{\mu } \right]}^{min\left\{ \left[ \frac{N-k}{\mu } \right],n \right\}}{}{{(-1)}^{l}}C_{n}^{l}{{b}_{\mu ,N}}(l\mu +k),\quad k=-1,-2,\ldots ,-\mu n,$$
$${{b}_{\mu ,N}}(k)=\sum\limits_{m=k}^{N}{}a(m){{d}_{\mu }}(m-k)=(D_{N}^{\mu }{{a}_{N}}{{)}_{k}},\quad k=0,1,\ldots ,N.$$
Here $D_{N}^{\mu }$ is the matrix of dimension $(N+1)\times (N+1)$ with elements $D_{k,j}^{\mu }={{d}_{\mu }}(j-k)$ if $0\le k\le j\le N$, and $D_{k,j}^{\mu }=0$ if $j<k$ or $j,k>N$; ${{a}_{N}}=(a(0),a(1),a(2),\ldots ,a(N))$. The spectral characteristic of the optimal estimate ${{\hat{A}}_{N}}\xi $ can be calculated by the following formulas:
\begin{equation}\label{eq31}
h_{\mu ,N}^{(a)}(\lambda )={{A}_{N}}({{e}^{i\lambda }})-r_{\mu ,N}^{(a)}({{e}^{i\lambda }})\Phi _{\mu }^{-1}({{e}^{-i\lambda }}),
\end{equation}
\begin{equation}\label{eq32}
{{A}_{N}}({{e}^{i\lambda }})=\sum\limits_{k=0}^{N}{}{{(D_{N}^{\mu }{{a}_{N}})}_{k}}{{e}^{i\lambda k}},\quad r_{\mu ,N}^{(a)}({{e}^{i\lambda }})=\sum\limits_{j=0}^{N}{}{{(D_{N}^{\mu }{{A}_{N}}{{\varphi }_{\mu ,N}})}_{j}}{{e}^{i\lambda j}},
\end{equation}
where the matrix ${{A}_{N}}$ of dimension $(N+1)\times (N+1)$ is determined by coefficients $a(k)$, $k=0,1,\ldots ,N$, in the following way: ${{({{A}_{N}})}_{k,j}}=a(k+j)$ if $0\le k+j\le N$, ${{({{A}_{N}})}_{k,j}}=0$ if $k+j>N$, $0\le k,j\le N$. The value of the mean square error of the optimal estimate ${{\hat{A}}_{N}}\xi $ can be calculated by the following formula:
\begin{equation}\label{eq33}
\Delta (f,{{\hat{A}}_{N}}\xi ):=\text{E}|{{A}_{N}}\xi -{{\hat{A}}_{N}}\xi {{|}^{2}}=\frac{1}{2\pi }\int\limits_{-\pi }^{\pi }{|r_{\mu ,N}^{(a)}({{e}^{i\lambda }}{{)|}^{2}}d\lambda }=||D_{N}^{\mu }{{A}_{N}}{{\varphi }_{\mu ,N}}{{||}^{2}}.
\end{equation}

Consequently, the following theorem holds true.

\begin{thm} \label{thm3.3}
 Let a stochastic sequence $\{\xi (m),m\in \mathbb{Z}\}$ determine a stationary stochastic $n$th increment sequence ${{\xi }^{(n)}}(m,\mu )$ with absolutely continuous spectral function $F(\lambda )$ and spectral density $f(\lambda )$ which admits the canonical factorization $(12)$. The optimal linear estimate ${{\hat{A}}_{N}}\xi $ of the functional ${{A}_{N}}\xi $ of unobserved values $\xi (m)$, $m=0,1,2,\ldots $, from observations of the sequence $\xi (m)$ at points $m=-1,-2,\ldots $ can be calculated by formula$(30)$. The spectral characteristic $h_{\mu ,N}^{(a)}(\lambda )$ of the optimal linear estimate ${{\hat{A}}_{N}}\xi $ can be calculated by formula $(31)$. The value of mean square error $\Delta (f,{{\hat{A}}_{N}}\xi )$ can be calculated by formula $(33)$.
\end{thm}

Consider the case where $\mu >m\ge 0$. In this case the mean square optimal estimate of the value $\xi (m)$, $m\ge 0$, can be calculated by formula
\begin{equation}\label{eq34}
\hat{\xi }(m)=-\sum\limits_{l=1}^{n}{}{{(-1)}^{l}}C_{n}^{l}\xi (m-l\mu )+\int\limits_{-\pi }^{\pi }{{{h}_{m,\mu }}(\lambda )(1-{{e}^{-i\lambda \mu }}{{)}^{n}}\frac{1}{{{(i\lambda )}^{n}}}dZ(\lambda )}
\end{equation}
The spectral characteristic ${{h}_{m,\mu }}(\lambda )$ and the value of the mean square error $\Delta (f,\hat{\xi }(m))=\Delta (f,{{\hat{\xi }}^{(n)}}(m,\mu ))$ of the estimate of the element $\xi (m)$ can be calculated by formulas $(24)$ and $(25)$ respectively.

Consequently, the following statement holds true.

\begin{nas} \label{corr3.2}
Let $\mu >m\ge 0$.
The optimal mean square estimate $\hat{\xi }(m)$ of the element $\xi (m)$, $\mu >m\ge 0$, from observations of the sequence $\xi (m)$ at points $m=-1,-2,\ldots $ can be calculated by formula $(34)$. The spectral characteristic ${{h}_{m,\mu }}(\lambda )$ of the optimal linear estimate $\hat{\xi }(m)$ can be calculated by formula $(24)$. The value of mean square error $\Delta (f,\hat{\xi }(m))$ can be calculated by formula $(25)$.
\end{nas}

\begin{zau} \label{rem3.1}
 Using relation $(13)$ we can find a relationship between coefficients $\{{{\varphi }_{\mu }}(k):k=0,1,2,\ldots \}$ and $\{\varphi (k):k=0,1,2,\ldots \}$. So far as
$$\int\limits_{-\pi }^{\pi }{\left| ln\frac{|1-{{e}^{-i\lambda \mu }}{{|}^{2n}}}{{{\lambda }^{2n}}} \right|d\lambda }<\infty $$
for every $n\ge 1$ and $\mu \ge 1$, there is a function ${{w}_{\mu }}(z)=\sum\nolimits_{k=0}^{\infty }{}{{w}_{\mu }}(k){{z}^{k}}$ such that
$$\sum\nolimits_{k=0}^{\infty }{}|{{w}_{\mu }}(k{{)|}^{2}}<\infty,\quad \frac{|1-{{e}^{-i\lambda \mu }}{{|}^{2n}}}{{{\lambda }^{2n}}}=|{{w}_{\mu }}({{e}^{-i\lambda }}{{)|}^{2}}$$
 and the following representation holds true:
\begin{equation}\label{eq35}
{{\Phi }_{\mu }}({{e}^{-i\lambda }})={{w}_{\mu }}({{e}^{-i\lambda }})\Phi ({{e}^{-i\lambda }}).
\end{equation}
 The function ${{w}_{\mu }}(z)$ is determined by the relation
\begin{equation}\label{eq36}
{{w}_{\mu }}(z)=Exp\left\{ \frac{1}{4\pi }\int\limits_{-\pi }^{\pi }{\frac{{{e}^{i\lambda }}+z}{{{e}^{i\lambda }}-z}ln\frac{|1-{{e}^{-i\lambda \mu }}{{|}^{2n}}}{{{\lambda }^{2n}}}d\lambda } \right\}.
\end{equation}
 From $(35)$ we can get
$${{\varphi }_{\mu }}(k)=\sum\limits_{j=0}^{k}{}{{w}_{\mu }}(k-j)\varphi (j).$$
Therefore elements ${{\varphi }_{\mu }}=({{\varphi }_{\mu }}(0),{{\varphi }_{\mu }}(1),{{\varphi }_{\mu }}(2),\ldots )$ and $\varphi =(\varphi (0),\varphi (1),\varphi (2),\ldots )$ from the space ${{\ell }_{2}}$ are connected by the following relation
\begin{equation}\label{eq37}
{{\varphi }_{\mu }}={{W}^{\mu }}\varphi,
\end{equation}
where ${{W}^{\mu }}$ is a linear operator in the space ${{\ell }_{2}}$ with elements $W_{j,k}^{\mu }={{w}_{\mu }}(j-k)$ if $0\le k\le j$ and $W_{j,k}^{\mu }=0$ if $j<k.$ The vectors ${{\varphi }_{\mu ,N}}=({{\varphi }_{\mu }}(0),{{\varphi }_{\mu }}(1),{{\varphi }_{\mu }}(2),\ldots ,{{\varphi }_{\mu }}(N))$ and ${{\varphi }_{N}}=(\varphi (0),\varphi (1),\varphi (2),\ldots ,\varphi (N))$ are connected by the relation
\begin{equation}\label{eq38}
{{\varphi }_{\mu ,N}}=W_{N}^{\mu }{{\varphi }_{N}},
\end{equation}
where $W_{N}^{\mu }$ ia a matrix of dimension $(N+1)\times (N+1)$ with elements $W_{j,k}^{\mu }={{w}_{\mu }}(j-k)$ if $0\le k\le j\le N$ and $W_{j,k}^{\mu }=0$ if $j<k$, $j,k=0,1,\ldots ,N$.
\end{zau}

\begin{exm} \label{exm3.1}
 Consider an $ARIMA(0,1,1)$ sequence $\{\xi (m):m\in \mathbb{Z}\}$. Increments of order 1 of the sequence $\xi (m)$ are stationary and increments with step 1 form one-sided moving average stochastic sequence of order 1 with parameter $\phi $. The spectral density of the sequence $\xi (m)$can be expressed as
$$f(\lambda )=\frac{{{\lambda }^{2}}|1+\phi {{e}^{i\lambda }}{{|}^{2}}}{|1-{{e}^{i\lambda }}{{|}^{2}}}.$$
By using $(12)$ and $(13)$ the function ${{\Phi }_{\mu }}(\lambda )$, $\mu >1$, is calculated by formula
$${{\Phi }_{\mu }}(\lambda )=1+(1+\phi ){{e}^{-i\lambda }}+\ldots +(1+\phi ){{e}^{-i\lambda (\mu -1)}}+\phi {{e}^{-i\lambda \mu }}.$$
Thus increments of order 1 with step $\mu >0$ of the sequence $\xi (m)$ form one-sided moving average stochastic sequence of order $\mu $.

Consider the problem of finding the mean square optimal linear estimate of the value of the functional ${{A}_{1}}\xi =a\xi (0)+b\xi (1)$ which depends of unknown values $\xi (0)$, $\xi (1)$ of the stochastic sequence $\xi (m)$ from observations $\xi (m)$ at points $m=-1,-2,\ldots $. We use theorem $3.3$ to solve this problem. The spectral characteristic $(31)$ of the optimal estimate ${{\hat{A}}_{1}}\xi $ of the functional ${{A}_{1}}\xi $ can be calculated by the formula
$$h_{\mu ,1}^{(a)}(\lambda )=(a+{{\delta }_{\mu 1}}b)+b{{e}^{i\lambda }}-\frac{(1-{{e}^{-i\lambda }})(a+b(1+\phi )+b{{e}^{i\lambda }})}{(1+\phi {{e}^{-i\lambda }})(1-{{e}^{-i\lambda \mu }})},$$
where ${{\delta }_{\mu 1}}$ is the Kronecker symbol. Using formula $(30)$ we calculated an estimate of the functional ${{A}_{1}}\xi $
$${{\hat{A}}_{1}}\xi =(a+b)(1+\phi )\sum\limits_{k=1}^{\infty }{}{{(-\phi )}^{k-1}}\xi (-k).$$
The value of mean square error is calculated by formula $(33)$
$$\Delta (f,{{\hat{A}}_{1}}\xi )={{a}^{2}}+2ab(1+\phi )+{{b}^{2}}(2+2\phi +{{\phi }^{2}}).$$
\end{exm}

\section{{Minimax-robust method of extrapolation}}

The proposed formulas may be employed under the condition that the spectral density $f(\lambda )$ of the considered stochastic sequence $\xi (m)$ with stationary $n$th increments is known. The value of the mean square error $\Delta (h_{\mu }^{(a)}(f);f):=\Delta (f,\hat{A}\xi )$ and the spectral characteristic $h_{\mu }^{(a)}(f)$ of the optimal linear estimate $\hat{A}\xi $ of the functional $A\xi $ which depends of unknown values $\xi (m)$ can be calculated by formulas $(27)$ and $(29)$, the value of mean square error $\Delta (h_{\mu ,N}^{(a)}(f);f):=\Delta (f,{{\hat{A}}_{N}}\xi )$ and the spectral characteristic $h_{\mu ,N}^{(a)}(f)$ of the optimal linear estimate ${{\hat{A}}_{N}}\xi $ of the functional ${{A}_{N}}\xi $ which depends of unknown values $\xi (m)$ can be calculated by formulas $(31)$ and $(33)$. In the case where the spectral density is not exactly known, but a set $D$ of admissible spectral densities is given, the minimax (robust) approach to estimation of the functionals of the unknown values of a stochastic sequence with stationary increments is reasonable. In other words we are interesting in finding an estimate that minimizes the maximum of the mean square errors for all spectral densities from a given class $D$ of admissible spectral densities simultaneously.

\begin{ozn} \label{def4.1}
 For a given class of spectral densities $D$ a spectral density ${{f}_{0}}(\lambda )\in D$ is called least favorable in $D$ for the optimal linear estimate the functional $A\xi $ if the following relation holds true:
$$\Delta ({{f}_{0}})=\Delta (h_{\mu }^{(a)}({{f}_{0}});{{f}_{0}})=\underset{f\in D}{\max}\,\Delta (h_{\mu }^{(a)}(f);f).$$
\end{ozn}

\begin{ozn} \label{def4.2}
For a given class of spectral densities $D$ a spectral characteristic ${{h}^{0}}(\lambda )$ of the optimal linear estimate of the functional $A\xi $ is called minimax-robust if there are satisfied conditions
$${{h}^{0}}(\lambda )\in {{H}_{D}}=\bigcap\limits_{f\in D}{}L_{2}^{0-}(f),$$
$$\underset{h\in {{H}_{D}}}{\min}\,\underset{f\in D}{\max}\,\Delta (h;f)=\underset{f\in D}{\max }\,\Delta ({{h}^{0}};f).$$
\end{ozn}

Analyzing the derived formulas and using the introduced definitions we can conclude that the following statements are true.

\begin{lem} \label{lem4.1}
 Spectral density ${{f}^{0}}(\lambda )\in D$ which admits the canonical factorization $(12)$ is the least favorable in the class of admissible spectral densities $D$ for the optimal linear estimation of the functional $A\xi $ if
\begin{equation}\label{eq39}
{{f}^{0}}(\lambda )={{\left| \sum\limits_{k=0}^{\infty }{}{{\varphi }^{0}}(k){{e}^{-i\lambda k}} \right|}^{2}},
\end{equation}
where ${{\varphi }^{0}}=\{{{\varphi }^{0}}(k):k=0,1,2,\ldots \}$ is a solution to the conditional extremum problem
\begin{equation}\label{eq40}
||{{D}^{\mu }}A{{\varphi }_{\mu }}{{||}^{2}}\to \max,\quad f(\lambda )={{\left| \sum\limits_{k=0}^{\infty }{}\varphi (k){{e}^{-i\lambda k}} \right|}^{2}}\in D.
\end{equation}
\end{lem}

\begin{lem} \label{lem4.2}
 Spectral density ${{f}^{0}}(\lambda )\in D$ which admits the canonical factorization $(12)$ is the least favorable in the class of admissible spectral densities $D$ for the optimal linear estimation of the functional ${{A}_{N}}\xi $ if
\begin{equation}\label{eq41}
{{f}^{0}}(\lambda )={{\left| \sum\limits_{k=0}^{N}{}{{\varphi }^{0}}(k){{e}^{-i\lambda k}} \right|}^{2}},
\end{equation}
where $\varphi _{N}^{0}=\{{{\varphi }^{0}}(k):k=0,1,2,\ldots ,N\}$ is a solution to the conditional extremum problem
\begin{equation}\label{eq42}
||D_{N}^{\mu }{{A}_{N}}{{\varphi }_{\mu ,N}}{{||}^{2}}\to \max,\quad f(\lambda )={{\left| \sum\limits_{k=0}^{N}{}\varphi (k){{e}^{-i\lambda k}} \right|}^{2}}\in D.
\end{equation}
 If $h_{\mu }^{(a)}({{f}^{0}})\in {{H}_{D}}$, the minimax-robust spectral characteristic can be calculated as ${{h}^{0}}=h_{\mu }^{(a)}({{f}^{0}})$.
\end{lem}

The minimax-robust spectral characteristic ${{h}^{0}}$ and the least favorable spectral density ${{f}^{0}}$ form a saddle point of the function $\Delta (h;f)$ on the set ${{H}_{D}}\times D$. The saddle point inequalities
$$\Delta (h;{{f}^{0}})\ge \Delta ({{h}^{0}};{{f}^{0}})\ge \Delta ({{h}^{0}};f)\quad \forall f\in D\ \forall h\in {{H}_{D}}$$
hold true if ${{h}^{0}}=h_{\mu }^{(a)}({{f}^{0}})$ and $h_{\mu }^{(a)}({{f}^{0}})\in {{H}_{D}}$, where ${{f}^{0}}$ is a solution to the conditional extremum problem
\begin{equation}\label{eq43}
\tilde{\Delta }(f)=-\Delta (h_{\mu }^{(a)}({{f}^{0}});f)\to \inf,\quad f\in D,
\end{equation}
$$\Delta (h_{\mu }^{(a)}({{f}^{0}});f)=\frac{1}{2\pi }\int\limits_{-\pi }^{\pi }{\frac{|{{r}_{\mu }}({{e}^{i\lambda }}{{)|}^{2}}}{{{f}^{0}}(\lambda )}f(\lambda )d\lambda },$$
where ${{r}_{\mu }}({{e}^{i\lambda }})$ is determined by formula $(28)$ or $(32)$ with $f(\lambda )={{f}^{0}}(\lambda )$.
The conditional extremum problem (43) is equivalent to the unconditional extremum problem
	$${{\Delta }_{D}}(f)=\tilde{\Delta }(f)+\delta (f|D)\to \inf,$$
where $\delta (f|D)$ is the indicator function of the set $D$. Solution ${{f}^{0}}$ to this unconditional extremum problem is characterized by the condition $0\in \partial {{\Delta }_{D}}({{f}^{0}})$, where $\partial {{\Delta }_{D}}({{f}^{0}})$ is the subdifferential of the functional ${{\Delta }_{D}}({{f}^{0}})$ at point ${{f}^{0}}$ (see Pshenichnyi (1982) or Moklyachuk (2008)). With the help of the condition $0\in \partial {{\Delta }_{D}}({{f}^{0}})$ we can find the least favorable spectral densities in some special classes of spectral densities (see books by Moklyachuk (2008), Moklyachuk and Masyutka (2012) for more details).

\section{{Least favorable spectral densities in the class ${{D}_{0}}$}}

Consider the problem of the optimal estimation of functionals $A\xi $ and ${{A}_{N}}\xi $ of unknown values $\xi (k)$, $k=0,1,2\ldots $, of the stochastic sequence $\xi (k)$ with stationary $n$th increments in the case where the spectral density is not known, but the following set of spectral densities is given
$${{D}_{0}}=\left\{ f(\lambda )|\frac{1}{2\pi }\int\limits_{-\pi }^{\pi }{f(\lambda )d\lambda \le {{P}_{0}}} \right\}.$$
It follows from the condition $0\in \partial {{\Delta }_{D}}({{f}_{0}})$ for $D={{D}_{0}}$ that the least favorable density satisfies the equation
$$|r_{\mu }^{(a)}({{e}^{i\lambda }}{{)|}^{2}}{{({{f}^{0}}(\lambda ))}^{-1}}=\psi (\lambda )+{{c}^{-2}},$$
 where $\psi (\lambda )\le 0$ and $\psi (\lambda )=0$ if ${{f}^{0}}(\lambda )>0$. Therefore, the least favorable density in the class ${{D}_{0}}$ for the optimal linear estimation of the functional $A\xi $ can be presented in the form
\begin{equation}\label{eq44}
{{f}^{0}}(\lambda )={{\left| c\sum\limits_{k=0}^{\infty }{}{{({{D}^{\mu }}A{{\varphi }_{\mu }})}_{k}}{{e}^{i\lambda k}} \right|}^{2}},
\end{equation}
where the unknown parameters $c$, ${{\varphi }_{\mu }}=({{\varphi }_{\mu }}(0),{{\varphi }_{\mu }}(1),{{\varphi }_{\mu }}(2),\ldots )$ can be calculated using factorization $(12)$, equation $(37)$, condition $(40)$ and condition $\int_{-\pi }^{\pi }{}f(\lambda )d\lambda =2\pi {{P}_{0}}.$
Consider the equation
\begin{equation}\label{eq45}
{{D}^{\mu }}A{{W}^{\mu }}\varphi =\alpha \overline{\varphi },\quad \alpha \in C.
\end{equation}
For each solution of this equation such that $||\varphi {{||}^{2}}={{P}_{0}}$ the following equality holds true:
$${{f}^{0}}(\lambda )={{\left| \sum\limits_{k=0}^{\infty }{}\varphi (k){{e}^{i\lambda k}} \right|}^{2}}={{\left| c\sum\limits_{k=0}^{\infty }{}{{({{D}^{\mu }}A{{W}^{\mu }}\varphi )}_{k}}{{e}^{i\lambda k}} \right|}^{2}}.$$

Denote by ${{\nu }_{0}}{{P}_{0}}$ the maximum value of $||{{D}^{\mu }}A{{W}^{\mu }}\varphi {{||}^{2}}$ on the set of those solutions $\varphi $ of equation $(45)$, which satisfy condition $||\varphi {{||}^{2}}={{P}_{0}}$ and define canonical factorization $(12)$ of the spectral density ${{f}^{0}}(\lambda )$. Let $\nu _{0}^{+}{{P}_{0}}$ be the maximum value of $||{{D}^{\mu }}A{{W}^{\mu }}\varphi {{||}^{2}}$ on the set of those $\varphi $ which satisfy condition $||\varphi {{||}^{2}}={{P}_{0}}$ and define canonical factorization $(12)$ of the spectral density ${{f}^{0}}(\lambda )$ defined by $(44)$.

The derived equations and conditions give us a possibility to verify the validity of following statement.

\begin{thm} \label{thm5.1}
If there exists a solution ${{\varphi }^{0}}=\{{{\varphi }^{0}}(m):m\ge 0\}$ of equation $(45)$ which satisfies conditions $||{{\varphi }^{0}}{{||}^{2}}={{P}_{0}}$ and ${{\nu }_{0}}{{P}_{0}}=\nu _{0}^{+}{{P}_{0}}=||{{D}^{\mu }}A{{W}^{\mu }}{{\varphi }^{0}}{{||}^{2}}$, the spectral density $(39)$ is least favorable density in the class ${{D}_{0}}$ for the optimal estimation of the functional $A\xi $ of unknown values $\xi (k)$, $k=0,1,2\ldots $, of the stochastic sequence $\xi (m)$ with stationary $n$th increments. The increment sequence  ${{\xi }^{(n)}}(m,\mu )$ admits a one-sided moving average representation. If ${{\nu }_{0}}<\nu _{0}^{+}$, the density $(44)$ which admits the canonical factorization $(12)$ is least favorable in the class ${{D}_{0}}$. The sequence $c{{\varphi }_{\mu }}=\{c{{\varphi }_{\mu }}(k):k\ge 0\}$ is determined by equality $(37)$, conditions $(40)$ and the condition $\int_{-\pi }^{\pi }{}f(\lambda )d\lambda =2\pi {{P}_{0}}$.
\end{thm}

Consider the problem of optimal estimation of the functional ${{A}_{N}}\xi $. In this case the least favorable spectral density is determined by the relation
\begin{equation}\label{eq46}
{{f}^{0}}(\lambda )={{\left| c\sum\limits_{k=0}^{N}{}{{(D_{N}^{\mu }{{A}_{N}}{{\varphi }_{\mu ,N}})}_{k}}{{e}^{i\lambda k}} \right|}^{2}}.
\end{equation}
Define the matrix $\hat{D}_{N}^{\mu }$ with the help of relation
\begin{equation}\label{eq47}
{{(\hat{D}_{N}^{\mu }{{A}_{N}}{{\varphi }_{\mu ,N}})}_{k}}=\sum\limits_{m=0}^{N}{}\sum\limits_{l=N-k}^{N}{}{{\varphi }_{\mu }}(m)a(m+l){{d}_{\mu }}(l+k-N),\quad k=0,1,2,\ldots ,N,
\end{equation}
where $a(p)=0$ if $p>N$. Taking into consideration $(38)$, we have the following equality
\begin{equation}\label{eq48}
{{\left| r_{\mu ,N}^{(a)}({{e}^{i\lambda }}) \right|}^{2}}={{\left| \sum\limits_{j=0}^{N}{}{{(D_{N}^{\mu }{{A}_{N}}W_{N}^{\mu }{{\varphi }_{N}})}_{j}}{{e}^{i\lambda j}} \right|}^{2}}={{\left| \sum\limits_{j=0}^{N}{}{{(\hat{D}_{N}^{\mu }{{A}_{N}}W_{N}^{\mu }{{\varphi }_{N}})}_{j}}{{e}^{-i\lambda j}} \right|}^{2}}.
\end{equation}
Therefore each solution ${{\varphi }_{N}}=({{\varphi }^{0}}(0),{{\varphi }^{0}}(1),{{\varphi }^{0}}(2),\ldots ,{{\varphi }^{0}}(N))$ of the equation
\begin{equation}\label{eq49}
D_{N}^{\mu }{{A}_{N}}W_{N}^{\mu }{{\varphi }_{N}}=\alpha {{\overline{\varphi }}_{N}},\quad \alpha \in C,
\end{equation}
 or the equation
 \begin{equation}\label{eq50}
\hat{D}_{N}^{\mu }{{A}_{N}}W_{N}^{\mu }{{\varphi }_{N}}=\beta {{\overline{\varphi }}_{N}},\quad \beta \in C,
\end{equation}
 such that $||{{\varphi }_{N}}{{||}^{2}}={{P}_{0}}$, satisfies the following equality
\begin{equation}\label{eq51}
{{f}^{0}}(\lambda )={{\left| \sum\limits_{k=0}^{N}{}\varphi (k){{e}^{i\lambda k}} \right|}^{2}}={{\left| r_{\mu ,N}^{(a)}({{e}^{i\lambda }}) \right|}^{2}}.
\end{equation}

Denote by $\nu _{0}^{N}{{P}_{0}}$ the maximum value of $||D_{N}^{\mu }{{A}_{N}}W_{N}^{\mu }{{\varphi }_{N}}{{||}^{2}}=||\hat{D}_{N}^{\mu }{{A}_{N}}W_{N}^{\mu }{{\varphi }_{N}}{{||}^{2}}$ on the set of solutions ${{\varphi }_{N}}$ of equation $(49)$ or equation $(50)$, which satisfy condition $||{{\varphi }_{N}}{{||}^{2}}={{P}_{0}}$ and determine the canonical factorization $(12)$ of the spectral density ${{f}^{0}}(\lambda )\in {{D}_{0}}$. Let $\nu _{0}^{N+}{{P}_{0}}$ be the maximum value of $||D_{N}^{\mu }{{A}_{N}}W_{N}^{\mu }{{\varphi }_{N}}{{||}^{2}}$ on the set of those ${{\varphi }_{N}}$ which satisfy condition $||{{\varphi }_{N}}{{||}^{2}}={{P}_{0}}$ and determine the canonical factorization ($12$) of the spectral density ${{f}^{0}}(\lambda )$ defined by ($46$).
The following statement holds true.

\begin{thm} \label{thm5.2}
 If there exists a solution $\varphi _{N}^{0}=\{{{\varphi }^{0}}(m):m=0,1,2,\ldots ,N\}$ of equation $(49)$ or equation $(50)$ which satisfies conditions $||\varphi _{N}^{0}{{||}^{2}}={{P}_{0}}$ and ${{\nu }_{0}}{{P}_{0}}=\nu _{0}^{+}{{P}_{0}}=||D_{N}^{\mu }{{A}_{N}}W_{N}^{\mu }\varphi _{N}^{0}{{||}^{2}}$, the spectral density $(41)$ is least favorable in the class ${{D}_{0}}$ for the optimal estimation of the functional ${{A}_{N}}\xi $ of unknown values $\xi (k)$, $k=0,1,2\ldots ,N$, of the stochastic sequence $\xi (m)$ with stationary $n$th increments. The increment sequence  ${{\xi }^{(n)}}(m,\mu )$ admits a one-sided moving average representation of order $N$. If ${{\nu }_{0}}<\nu _{0}^{+}$, the density $(46)$ which admits the canonical factorization $(12)$ is least favorable in the class ${{D}_{0}}$. The sequence $c{{\varphi }_{\mu ,N}}=\{c{{\varphi }_{\mu }}(k):k=0,1,2,\ldots ,N\}$ is determined by equation $(38)$, conditions $(42)$ and the condition $\int_{-\pi }^{\pi }{}f(\lambda )d\lambda =2\pi {{P}_{0}}$.
\end{thm}

\begin{exm} \label{exm5.1}  Consider the problem of minimax estimation of the functional ${{A}_{1}}\xi =a\xi (0)+b\xi (1)$ of a stochastic sequence $\{\xi (m):m\in \mathbb{Z}\}$ with stationary increments of order $1$ from observations of the sequence $\xi (m)$ for $m=-1,-2,\ldots $. We use theorem $5.2$ to solve this problem. The matrices used in $(42)$ and $(38)$ are the following
$${{A}_{1}}=\left( \begin{matrix}
   a & b  \\
   b & 0  \\
\end{matrix} \right),\quad D_{1}^{\mu }=\left( \begin{matrix}
   1 & {{\delta }_{\mu 1}}  \\
   0 & 1  \\
\end{matrix} \right),\quad W_{1}^{\mu }=\left( \begin{matrix}
   {{w}_{\mu }}(0) & 0  \\
   {{w}_{\mu }}(1) & {{w}_{\mu }}(0)  \\
\end{matrix} \right),
$$
where ${{w}_{\mu }}(0)$, ${{w}_{\mu }}(1)$ are the Fourier coefficients of the function ${{w}_{\mu }}({{e}^{-i\lambda }})$ defined by $(36)$. The least favorable density in the set ${{D}_{0}}$ is defined by a solution of the optimization problem $(42)$, where ${{\varphi }_{\mu ,N}}={{\varphi }_{\mu ,1}}=D_{1}^{\mu }{{\varphi }_{1}}$, ${{\varphi }_{1}}=({{\varphi }_{1}}(0),{{\varphi }_{1}}(1){)}'$. Let us assume that $xy\ne 0$, where $x:=(a+{{\delta }_{\mu 1}}b){{w}_{\mu }}(0)+b{{w}_{\mu }}(1)$, $y:=b{{w}_{\mu }}(0)$. Then the optimization problem can be represented in the form
$$\left\{ \begin{matrix}
   {{(x\varphi (0)+y\varphi (1))}^{2}}+{{y}^{2}}{{\varphi }^{2}}(0)\to max;  \\
   {{\varphi }^{2}}(0)+{{\varphi }^{2}}(1)\le {{P}_{0}},  \\
\end{matrix} \right.$$
A solution $\varphi _{1}^{0}=({{\varphi }^{0}}(0),{{\varphi }^{0}}(1){)}'$ of this problem is calculated as follows
$${{\varphi }^{0}}(0)=\pm {{\left( \frac{{{P}_{0}}({{x}^{2}}+4{{y}^{2}}+\sqrt{{{x}^{2}}({{x}^{2}}+4{{y}^{2}})})}{2({{x}^{2}}+4{{y}^{2}})} \right)}^{\frac{1}{2}}};$$
$${{\varphi }^{0}}(1)=\pm sign(xy){{\left( \frac{{{P}_{0}}({{x}^{2}}+4{{y}^{2}}-\sqrt{{{x}^{2}}({{x}^{2}}+4{{y}^{2}})})}{2({{x}^{2}}+4{{y}^{2}})} \right)}^{\frac{1}{2}}}.$$
The vector $\varphi _{1}^{0}=({{\varphi }^{0}}(0),{{\varphi }^{0}}(1){)}'$ provides the maximum value of $||D_{1}^{\mu }{{A}_{1}}W_{1}^{\mu }{{\varphi }_{1}}{{||}^{2}}$, satisfies condition $||\varphi _{1}^{0}{{||}^{2}}={{P}_{0}}$ and equation $(49)$ with $\lambda =\frac{x+\sqrt{{{x}^{2}}+4{{y}^{2}}}}{2y}$ if $y>0$, and with $\lambda =\frac{x-\sqrt{{{x}^{2}}+4{{y}^{2}}}}{2y}$ if $y<0$. Using theorem $5.2$ we can conclude that the spectral density ${{f}^{0}}(\lambda )=|{{\varphi }^{0}}(0)+{{\varphi }^{0}}(1){{e}^{-i\lambda }}{{|}^{2}}$ is the least favorable one in the class ${{D}_{0}}$ for the optimal estimation of the functional ${{A}_{1}}\xi =a\xi (0)+b\xi (1)$ of unknown values $\xi (0)$, $\xi (1)$ of the stochastic sequence $\xi (m)$ with stationary $n$th increments.
\end{exm}

\section{{Least favorable spectral densities in the class ${{D}_{M}}$}}

Consider the problem of the optimal estimation of functionals $A\xi $ and ${{A}_{N}}\xi $ of unknown values $\xi (k)$, $k=0,1,2\ldots $, of the stochastic sequence $\xi (k)$ with stationary $n$th increments in the case where the spectral density is not exactly known, but the following set of spectral densities is given
$${{D}_{M}}=\left\{ f(\lambda )|\frac{1}{2\pi }\int\limits_{-\pi }^{\pi }{f(\lambda )cos(m\lambda )d\lambda }={{\rho }_{m}},m=0,1,2,\ldots ,M \right\},$$
where ${{\rho }_{0}}={{P}_{0}}$ and $\{{{\rho }_{m}}$, $m=0,1,2,\ldots ,M\}$ is a strictly positive sequence (see Krein and Nudel'man (1977)). It follows from the condition $0\in \partial {{\Delta }_{D}}({{f}_{0}})$ that the least favorable density satisfies the equation
$$|r_{\mu }^{(a)}({{e}^{i\lambda }}{{)|}^{2}}{{({{f}^{0}}(\lambda ))}^{-1}}=\psi (\lambda )+c\sum\limits_{m=1}^{M}{}{{\psi }_{m}}cosm\lambda .$$
Thus, the least favorable density in the class ${{D}_{M}}$ for the optimal linear estimation of the functional $A\xi $ can be presented in the form
\begin{equation}\label{eq52}
{{f}^{0}}(\lambda )=\frac{{{\left| {{c}_{0}}\sum\limits_{k=0}^{\infty }{}{{({{D}^{\mu }}A{{\varphi }_{\mu }})}_{k}}{{e}^{i\lambda k}} \right|}^{2}}}{{{\left| \sum\limits_{k=1}^{M}{}{{c}_{m}}{{e}^{-i\lambda k}} \right|}^{2}}},
\end{equation}
where parameters ${{c}_{m}}$, $m=0,1,2,\ldots ,M$, ${{\varphi }_{\mu }}=({{\varphi }_{\mu }}(0),{{\varphi }_{\mu }}(1),{{\varphi }_{\mu }}(2),\ldots )$ can be calculated using conditions $(40)$, condition $\int_{-\pi }^{\pi }{}f(\lambda )cos(m\lambda )d\lambda =2\pi {{\rho }_{m}}$, $m=0,1,2,\ldots ,M$, equation $(37)$, factorization $(12)$.

Denote by ${{\nu }_{M}}{{P}_{0}}$ the maximum value of $||{{D}^{\mu }}A{{W}^{\mu }}\varphi {{||}^{2}}$ on the set of solutions $\varphi $ of the equation $(45)$ which satisfy condition $||\varphi {{||}^{2}}={{P}_{0}}$ and determine the canonical factorization ($12$) of the spectral density ${{f}^{0}}(\lambda )$. Let $\nu _{M}^{+}{{P}_{0}}$ be the maximum value of $||{{D}^{\mu }}A{{W}^{\mu }}\varphi {{||}^{2}}$ on the set of those $\varphi $ which satisfy condition $||\varphi {{||}^{2}}={{P}_{0}}$ and determine the canonical factorization ($12$) of the spectral density ${{f}^{0}}(\lambda )\in {{D}_{M}}$ defined by $(52)$. The derived equations and conditions give us a possibility to verify the validity of following statement.

\begin{thm} \label{thm6.1}
 If there exists a solution ${{\varphi }^{0}}=\{{{\varphi }^{0}}(m):m\ge 0\}$ of equation $(45)$ which satisfies conditions $||{{\varphi }^{0}}{{||}^{2}}={{P}_{0}}$ and ${{\nu }_{0}}{{P}_{0}}=\nu _{M}^{+}{{P}_{0}}=||{{D}^{\mu }}A{{W}^{\mu }}{{\varphi }^{0}}{{||}^{2}}$, the spectral density $(39)$ is least favorable in the class ${{D}_{M}}$ for the optimal extrapolation of the functional $A\xi $ of unknown values$\xi (k)$, $k=0,1,2\ldots $, of the stochastic sequence with stationary $n$th increments. If ${{\nu }_{M}}<\nu _{M}^{+}$, the density $(52)$ which admits the canonical factorization $(12)$ is least favorable in the class ${{D}_{M}}$. The sequence ${{\varphi }_{\mu }}=\{{{\varphi }_{\mu }}(k):k\ge 0\}$ and unknown parameters ${{c}_{m}}$, $m=0,1,2,\ldots ,M$, are determined by equation $(37)$, conditions $(40)$ and conditions $\int_{-\pi }^{\pi }{}f(\lambda )cos(m\lambda )d\lambda =2\pi {{\rho }_{m}}$, $m=0,1,2,\ldots ,M$.
\end{thm}

In the case of estimation of the functional ${{A}_{N}}\xi $ the least favorable spectral density is defined by equation
\begin{equation}\label{eq53}
{{f}^{0}}(\lambda )=\frac{{{\left| {{c}_{0}}\sum\limits_{k=0}^{N}{}{{(D_{N}^{\mu }{{A}_{N}}{{\varphi }_{\mu ,N}})}_{k}}{{e}^{i\lambda k}} \right|}^{2}}}{{{\left| \sum\limits_{k=1}^{M}{}{{c}_{m}}{{e}^{-i\lambda k}} \right|}^{2}}}.
\end{equation}
Let the matrix $\hat{D}_{N}^{\mu }$ be defined by equality $(47)$. Then equality$(48)$ holds true. Therefore each solution ${{\varphi }_{N}}=({{\varphi }^{0}}(0),{{\varphi }^{0}}(1),{{\varphi }^{0}}(2),\ldots ,{{\varphi }^{0}}(N))$ of the equation $(49)$ or the equation $(50)$ such that $||{{\varphi }_{N}}{{||}^{2}}={{P}_{0}}$ satisfies equality $(51)$.

Denote by $\nu _{0}^{N}{{P}_{0}}$ be the maximum value of $||D_{N}^{\mu }{{A}_{N}}W_{N}^{\mu }{{\varphi }_{N}}{{||}^{2}}=||\hat{D}_{N}^{\mu }{{A}_{N}}W_{N}^{\mu }{{\varphi }_{N}}{{||}^{2}}$ on the set of solutions ${{\varphi }_{N}}$ of equation $(49)$ or equation $(50)$, which satisfy condition $||{{\varphi }_{N}}{{||}^{2}}={{P}_{0}}$ and determine the canonical factorization $(12)$ of the spectral density ${{f}^{0}}(\lambda )$. Let $\nu _{0}^{N+}{{P}_{0}}$ be the maximum value of $||D_{N}^{\mu }{{A}_{N}}W_{N}^{\mu }{{\varphi }_{N}}{{||}^{2}}$ on the set of those ${{\varphi }_{N}}$ which satisfy condition $||{{\varphi }_{N}}{{||}^{2}}={{P}_{0}}$ and determine the canonical factorization $(12)$ of the spectral density ${{f}^{0}}(\lambda )$ defined by $(53)$.

The following statement holds true.

\begin{thm} \label{thm6.2}
If there exists a solution $\varphi _{N}^{0}=\{{{\varphi }^{0}}(m):m=0,1,2,\ldots ,N\}$ of equation $(49)$ or equation $(50)$ which satisfies conditions $||\varphi _{N}^{0}{{||}^{2}}={{P}_{0}}$ and ${{\nu }_{0}}{{P}_{0}}=\nu _{0}^{+}{{P}_{0}}=||D_{N}^{\mu }{{A}_{N}}W_{N}^{\mu }\varphi _{N}^{0}{{||}^{2}}$, the spectral density $(41)$ is least favorable in the class ${{D}_{M}}$ for the optimal estimation of the functional ${{A}_{N}}\xi $ of unknown values $\xi (k)$, $k=0,1,2\ldots ,N$, of the stochastic sequence with stationary $n$th increments. The increment ${{\xi }^{(n)}}(m,\mu )$ admits auto-regressive moving average representation of order $(M,N)$. If ${{\nu }_{0}}<\nu _{0}^{+}$, the density $(53)$ which admits the canonical factorization $(12)$ is the least favorable in the class ${{D}_{M}}$. The unknown parameters ${{\varphi }_{\mu ,N}}=\{{{\varphi }_{\mu }}(k):k=0,1,2,\ldots ,N\}$ and ${{c}_{m}}$, $m=0,1,2,\ldots ,M$, are determined by equality $(38)$, conditions $(42)$ and conditions $\int_{-\pi }^{\pi }{}f(\lambda )cos(m\lambda )d\lambda =2\pi {{\rho }_{m}}$, $m=0,1,2,\ldots ,M$.
\end{thm}

\section{{Least favorable spectral densities in the class $D_{v}^{u}$}}

Consider the problem of the optimal estimation of functionals $A\xi $ and ${{A}_{N}}\xi $ of unknown values $\xi (k)$, $k=0,1,2\ldots $, of the stochastic sequence $\xi (k)$ with stationary $n$th increments in the case where the spectral density is not known, but the following set of spectral densities is given
$$D_{v}^{u}=\left\{ f(\lambda )|v(\lambda )\le f(\lambda )\le u(\lambda ),\frac{1}{2\pi }\int\limits_{-\pi }^{\pi }{f(\lambda )d\lambda }\le {{P}_{0}} \right\},$$
 here $v(\lambda )$ and $u(\lambda )$ are some given (fixed) spectral densities. It follows from the condition $0\in \partial {{\Delta }_{D}}({{f}_{0}})$ for $D=D_{v}^{u}$ that the least favorable density ${{f}_{0}}$ in the class $D_{v}^{u}$ for the optimal linear estimation of the functional $A\xi $ is of the form
\begin{equation}\label{eq54}
{{f}^{0}}(\lambda )=max\left\{ v(\lambda ),min\left\{ u(\lambda ),{{\left| c\sum\limits_{k=0}^{\infty }{}{{({{D}^{\mu }}A{{\varphi }_{\mu }})}_{k}}{{e}^{i\lambda k}} \right|}^{2}} \right\} \right\},
\end{equation}
where the unknown parameters $c$, ${{\varphi }_{\mu }}=({{\varphi }_{\mu }}(0),{{\varphi }_{\mu }}(1),{{\varphi }_{\mu }}(2),\ldots )$ can be calculated using factorization $(12)$, equation $(37)$, conditions $(40)$ and condition $\int_{-\pi }^{\pi }{}f(\lambda )d\lambda =2\pi {{P}_{0}}.$

Denote by ${{\nu }_{u}}{{P}_{0}}$ the maximum value of $||{{D}^{\mu }}A{{W}^{\mu }}\varphi {{||}^{2}}$ on the set of those solutions $\varphi $ of equation $(45)$, which satisfy inequalities
$$v(\lambda )\le {{\left| \sum\limits_{k=0}^{\infty }{}\varphi (k){{e}^{-i\lambda k}} \right|}^{2}}\le u(\lambda ),$$
satisfy condition $||\varphi {{||}^{2}}={{P}_{0}}$ and determine the canonical factorization $(12)$ of the spectral density ${{f}^{0}}(\lambda )$. Let $\nu _{u}^{+}{{P}_{0}}$ be the maximum value of $||{{D}^{\mu }}A{{W}^{\mu }}\varphi {{||}^{2}}$ on the set of those $\varphi $ which satisfy condition $||\varphi {{||}^{2}}={{P}_{0}}$ and determine the canonical factorization $(12)$ of the spectral density ${{f}^{0}}(\lambda )$ defined by $(54)$. The derived equations and conditions give us a possibility to verify the validity of the following statement.

\begin{thm} \label{thm7.1}
If there exists a solution ${{\varphi }^{0}}=\{{{\varphi }^{0}}(m):m\ge 0\}$ of equation $(45)$ which satisfies conditions $||{{\varphi }^{0}}{{||}^{2}}={{P}_{0}}$ and ${{\nu }_{u}}{{P}_{0}}=\nu _{u}^{+}{{P}_{0}}=||{{D}^{\mu }}A{{W}^{\mu }}{{\varphi }^{0}}{{||}^{2}}$, the spectral density $(39)$ is least favorable in the class $D_{v}^{u}$ for the optimal estimation of the functional $A\xi $ of unknown values $\xi (k)$,$k=0,1,2\ldots $, of the stochastic sequence $\xi (m)$ with stationary $n$th increments. The increment sequence ${{\xi }^{(n)}}(m,\mu )$ admits one-sided moving average representation. If ${{\nu }_{u}}<\nu _{u}^{+}$, the density $(54)$ which admits the canonical factorization $(12)$ is least favorable in the class $D_{v}^{u}$. The sequence $c{{\varphi }_{\mu }}=\{c{{\varphi }_{\mu }}(k):k\ge 0\}$ is determined by equality $(37)$, conditions $(40)$ and the condition $\int_{-\pi }^{\pi }{}f(\lambda )d\lambda =2\pi {{P}_{0}}$. The minimax-robust spectral characteristic is calculated by formulas $(27)$, $(28)$.
\end{thm}

Consider the problem of the optimal estimation of the functional ${{A}_{N}}\xi $. In this case the least favorable spectral density is determined by the relation
\begin{equation}\label{eq54}
{{f}^{0}}(\lambda )=max\left\{ v(\lambda ),min\left\{ u(\lambda ),{{\left| c\sum\limits_{k=0}^{N}{}{{({{D}^{\mu }}A{{\varphi }_{\mu }})}_{k}}{{e}^{i\lambda k}} \right|}^{2}} \right\} \right\}.
\end{equation}
Denote by $\nu _{u}^{N}{{P}_{0}}$ the maximum value of $||D_{N}^{\mu }{{A}_{N}}W_{N}^{\mu }{{\varphi }_{N}}{{||}^{2}}=||\hat{D}_{N}^{\mu }{{A}_{N}}W_{N}^{\mu }{{\varphi }_{N}}{{||}^{2}}$ on the set of solutions ${{\varphi }_{N}}$ of equations $(49)$ and $(50)$ which satisfy inequality
$$v(\lambda )\le {{\left| \sum\limits_{k=0}^{N}{}\varphi (k){{e}^{-i\lambda k}} \right|}^{2}}\le u(\lambda ),
$$
satisfy condition $||{{\varphi }_{N}}{{||}^{2}}={{P}_{0}}$ and define the canonical factorization $(12)$ of the spectral density ${{f}^{0}}(\lambda )\in D_{v}^{u}$. Let $\nu _{u}^{N+}{{P}_{0}}$ be the maximum value of $||D_{N}^{\mu }{{A}_{N}}W_{N}^{\mu }{{\varphi }_{N}}{{||}^{2}}$ on the set of those ${{\varphi }_{N}}$ which satisfy condition $||{{\varphi }_{N}}{{||}^{2}}={{P}_{0}}$ and define canonical factorization $(12)$ of the spectral density ${{f}^{0}}(\lambda )$ determined by $(55)$.

The following statement holds true.

\begin{thm} \label{thm7.2}
 If there exists a solution $\varphi _{N}^{0}=\{{{\varphi }^{0}}(m):m=0,1,2,\ldots ,N\}$ of equation $(49)$ or equation $(50)$ which satisfies conditions $||\varphi _{N}^{0}{{||}^{2}}={{P}_{0}}$ and $\nu _{u}^{N}{{P}_{0}}=\nu _{u}^{N+}{{P}_{0}}=||D_{N}^{\mu }{{A}_{N}}W_{N}^{\mu }\varphi _{N}^{0}{{||}^{2}}$, spectral density $(41)$ is least favorable in the class $D_{v}^{u}$ for the optimal estimation of the functional ${{A}_{N}}\xi $ of unknown values $\xi (k)$, $k=0,1,2\ldots ,N$, of the stochastic sequence $\xi (m)$ with stationary $n$th increments. The increment ${{\xi }^{(n)}}(m,\mu )$ admits one-sided moving average representation of order $N$. If $\nu _{u}^{N}<\nu _{u}^{N+}$, the density $(55)$ which admits the canonical factorization $(12)$ is least favorable in the class $D_{v}^{u}$. The sequence $c{{\varphi }_{\mu ,N}}=\{c{{\varphi }_{\mu }}(k):k=0,1,2,\ldots ,N\}$ is determined by equation $(38)$, conditions $(42)$ and $\int_{-\pi }^{\pi }{}f(\lambda )d\lambda =2\pi {{P}_{0}}$. The minimax-robust spectral characteristic is calculated by formulas $(31)$, $(32)$.
\end{thm}

\begin{nas} \label{corr7.1}
Corollary 7.1 If we take $v(\lambda )=0$ and $u(\lambda )=\infty $, two previous theorems give us solutions to the problem of the minimax estimation of the functionals $A\xi $ and ${{A}_{N}}\xi $ for the set of spectral densities

	$${{D}_{0}}=\left\{ f(\lambda )|\frac{1}{2\pi }\int\limits_{-\pi }^{\pi }{f(\lambda )d\lambda }\le {{P}_{0}} \right\}.$$

\end{nas}

\section{{Least favorable spectral densities in the class ${{D}_{\varepsilon }}$}}

Consider the problem of the optimal estimation of functionals $A\xi $ and ${{A}_{N}}\xi $ of unknown values $\xi (k)$, $k=0,1,2\ldots $, of the stochastic sequence $\xi (k)$ with stationary $n$th increments in the case where the spectral density is not known, but the following set of spectral densities is given
$${{D}_{\varepsilon }}=\left\{ f(\lambda )|\frac{1}{2\pi }\int\limits_{-\pi }^{\pi }{|f(\lambda )-v(\lambda )|d\lambda }\le \varepsilon  \right\},$$
where $v(\lambda )$ is a bounded spectral density.

From the condition $0\in \partial {{\Delta }_{D}}({{f}_{0}})$ for $D={{D}_{\varepsilon }}$ we find the following equation to determine the least favorable spectral densities
\begin{equation}\label{eq55}
{{f}^{0}}(\lambda )=max\left\{ v(\lambda ),{{\left| c\sum\limits_{k=0}^{\infty }{}{{({{D}^{\mu }}A{{\varphi }_{\mu }})}_{k}}{{e}^{i\lambda k}} \right|}^{2}} \right\}.
\end{equation}
Let us define
\begin{equation}\label{eq57}
\frac{1}{2\pi }\int\limits_{-\pi }^{\pi }{f(\lambda )d\lambda }=\varepsilon +\frac{1}{2\pi }\int\limits_{-\pi }^{\pi }{v(\lambda )d\lambda }={{P}_{1}}.
\end{equation}

Let ${{\nu }_{\varepsilon }}{{P}_{1}}$ be the maximum value of $||{{D}^{\mu }}A{{W}^{\mu }}\varphi {{||}^{2}}$ on the set of those $\varphi $ which belongs to the set of solutions of equation $(45)$, satisfy the inequality
$$v(\lambda )\le {{\left| \sum\limits_{k=0}^{\infty }{}\varphi (k){{e}^{-i\lambda k}} \right|}^{2}},$$
satisfy condition $||\varphi {{||}^{2}}={{P}_{1}}$ and determine the canonical factorization $(12)$ of the spectral density ${{f}^{0}}(\lambda )$. Let $\nu _{\varepsilon }^{+}{{P}_{1}}$ be the maximum value of $||{{D}^{\mu }}A{{W}^{\mu }}\varphi {{||}^{2}}$ on the set of those $\varphi $ which satisfy condition $||\varphi {{||}^{2}}={{P}_{1}}$ and determine the canonical factorization $(12)$ of the spectral density ${{f}^{0}}(\lambda )$ defined by $(56)$. The following statement holds true.

\begin{thm} \label{thm8.1}
 If there exists a solution ${{\varphi }^{0}}=\{{{\varphi }^{0}}(m):m\ge 0\}$ of equation $(45)$ which satisfies conditions $||{{\varphi }^{0}}{{||}^{2}}={{P}_{1}}$ and ${{\nu }_{\varepsilon }}{{P}_{1}}=\nu _{\varepsilon }^{+}{{P}_{1}}=||{{D}^{\mu }}A{{W}^{\mu }}{{\varphi }^{0}}{{||}^{2}}$, the spectral density $(39)$ is least favorable in the class ${{D}_{\varepsilon }}$ for the optimal extrapolation of the functional $A\xi $ of unknown values $\xi (k)$, $k=0,1,2\ldots $, of the stochastic sequence $\xi (m)$ with stationary $n$th increments. The increment ${{\xi }^{(n)}}(m,\mu )$ admits one-sided moving average representation. If ${{\nu }_{u}}<\nu _{u}^{+}$, the density $(56)$ which admits the canonical factorization $(12)$ is least favorable in the class ${{D}_{\varepsilon }}$. The sequence $c{{\varphi }_{\mu }}=\{c{{\varphi }_{\mu }}(k):k\ge 0\}$ is determined by equality$(37)$, conditions $(40)$ and $\int_{-\pi }^{\pi }{}f(\lambda )d\lambda =2\pi \varepsilon +\int_{-\pi }^{\pi }{}v(\lambda )d\lambda $. The minimax-robust spectral characteristic is calculated by formulas $(27)$, $(28)$.
\end{thm}
In the case of optimal estimation of the functional ${{A}_{N}}\xi $ the least favorable spectral density is determined by formula
\begin{equation}\label{eq58}
{{f}^{0}}(\lambda )=\max\left\{ v(\lambda ),{{\left| c\sum\limits_{k=0}^{N}{}{{({{D}^{\mu }}A{{\varphi }_{\mu }})}_{k}}{{e}^{i\lambda k}} \right|}^{2}} \right\}.
\end{equation}

Let $\nu _{\varepsilon }^{N}{{P}_{1}}$ be the maximum value of $||D_{N}^{\mu }{{A}_{N}}W_{N}^{\mu }{{\varphi }_{N}}{{||}^{2}}=||\hat{D}_{N}^{\mu }{{A}_{N}}W_{N}^{\mu }{{\varphi }_{N}}{{||}^{2}}$ on the set of those ${{\varphi }_{N}}$ which belong to the set of solutions of equation $(49)$ or equation $(50)$, satisfy the inequality
$$v(\lambda )\le {{\left| \sum\limits_{k=0}^{N}{}\varphi (k){{e}^{-i\lambda k}} \right|}^{2}},$$
satisfy condition $||{{\varphi }_{N}}{{||}^{2}}={{P}_{1}}$ and determined the canonical factorization $(12)$ of the spectral density ${{f}^{0}}(\lambda )$, ${{f}^{0}}(\lambda )\in {{D}_{\varepsilon }}$. Let $\nu _{\varepsilon }^{N+}{{P}_{1}}$ be the maximum value of $||D_{N}^{\mu }{{A}_{N}}W_{N}^{\mu }{{\varphi }_{N}}{{||}^{2}}$ on the set of those ${{\varphi }_{N}}$ which satisfy condition $||{{\varphi }_{N}}{{||}^{2}}={{P}_{1}}$ and determined the canonical factorization $(12)$ of the spectral density ${{f}^{0}}(\lambda )$ defined by $(58)$. The following statement holds true.

\begin{thm} \label{thm8.2}
 If there exists a solution $\varphi _{N}^{0}=\{{{\varphi }^{0}}(m):m=0,1,2,\ldots ,N\}$ of equation $(49)$ or equation $(50)$ which satisfies conditions $||\varphi _{N}^{0}{{||}^{2}}={{P}_{1}}$ and ${{\nu }_{\varepsilon }}{{P}_{1}}=\nu _{\varepsilon }^{+}{{P}_{1}}=||{{D}^{\mu }}AW_{N}^{\mu }\varphi _{N}^{0}{{||}^{2}}$, the spectral density $(41)$ is least favorable in the class ${{D}_{\varepsilon }}$ for the optimal extrapolation of the functional ${{A}_{N}}\xi $ of unknown values $\xi (k)$, $k=0,1,2\ldots ,N$, of the stochastic sequence $\xi (m)$ with stationary $n$th increments. The increment ${{\xi }^{(n)}}(m,\mu )$ admits one-sided moving average representation of order $N$. If ${{\nu }_{\varepsilon }}<\nu _{\varepsilon }^{+}$, the density $(58)$ which admits the canonical factorization $(12)$ is least favorable in the class ${{D}_{\varepsilon }}$. The sequence $c{{\varphi }_{\mu ,N}}=\{c{{\varphi }_{\mu }}(k):k=0,1,2,\ldots ,N\}$ is determined by equation $(38)$, conditions $(42)$ and $\int_{-\pi }^{\pi }{}f(\lambda )d\lambda =2\pi \varepsilon +\int_{-\pi }^{\pi }{}v(\lambda )d\lambda $. The minimax-robust spectral characteristic is calculated by formulas $(31)$, $(32)$.
\end{thm}

\section{{Conclusions}}

In this article we describe methods of solution of the problem of optimal linear estimation of functionals which depend on unknown values of a stochastic sequence $\xi (m)$ with stationary $n$th increments. Estimates are based on observations of the sequence $\xi (t)$ at points $t=-1,-2,\ldots $. Formulas are derived for computing the value of the mean-square error and the spectral characteristic of the optimal linear estimate of functionals in the case of spectral certainty where the spectral density of the sequence is exactly known.

In the case of spectral uncertainty where the spectral density is not exactly known but, instead, a set of admissible spectral densities is specified, the minimax-robust method is applied. We propose a representation of the mean square error in the form of a linear functional in ${{L}_{1}}$ with respect to spectral densities, which allows us to solve the corresponding conditional extremum problem and describe the minimax (robust) estimates of the functional. Formulas that determine the least favorable spectral densities and minimax (robust) spectral characteristic of the optimal linear estimates of the functionals are derived  for some concrete classes of admissible spectral densities.


\begin{thebibliography}{99}

\bibitem{Dubov2012}
Dubovets'ka, I.I., Masyutka, O.Yu., Moklyachuk, M.P., 2012. Interpolation of periodically correlated stochastic sequences. Theory Probability and Mathematical Statistics 84, 43-156.

\bibitem{Dubov2013a}
Dubovets'ka, I.I., Moklyachuk, M.P., 2013. Extrapolation of  periodically correlated processes  from observations with noise. Theory Probability and Mathematical Statistics 88, 60-75.

\bibitem{Dubov2013b}
Dubovets'ka, I.I., Moklyachuk,M.P., 2013. Minimax estimation problem for periodically correlated stochastic processes. Journal of Mathematics and System Science 3(1), 26-30.

\bibitem{Grenander}
Grenander, U., 1957. A prediction problem in game theory. Arkiv f\"or Matematik 3, 371-379.

\bibitem{Franke}
Franke, J., 1985. Minimax robust prediction of discrete time series. Z. Wahrscheinlichkeitstheor. Verw. Gebiete 68, 337-364.

\bibitem{FrankePoor}
Franke, J., Poor, H.V., 1984. Minimax-robust filtering and finite-length robust predictors, In: Robust and Nonlinear Time Series Analysis. Lecture Notes in Statistics, Springer-Verlag 26, 87-126.

\bibitem{KassamPoor}
Kassam, S.A., Poor, H.V., 1985. Robust techniques for signal processing: A survey. Proceedings of the IEEE 73, 433-481.

\bibitem{Karhunen}
Karhunen, K., 1947. Uber lineare Methoden in der Wahrscheinlichkeitsrechnung. Annales Academiae Scientiarum Fennicae. Series A I. Mathematica 37, 3-79.

\bibitem{Kolmogorov}
Kolmogorov, A.N., 1992. Selected works of A. N. Kolmogorov. Vol. II: Probability theory and mathematical statistics. Ed. by A. N. Shiryayev. Mathematics and Its Applications. Soviet Series. 26. Dordrecht etc.: Kluwer Academic Publishers, 584.

\bibitem{Krein}
Krein, M.G., Nudel'man, A.A., 1977. The Markov moment problem and extremal problems. Ideas and problems of P. L. Chebysev and A. A. Markov and their further development. Translations of Mathematical Monographs. Vol. 50. Providence, R.I.: American Mathematical Society(AMS), 552.

\bibitem{Luz2012a}
Luz, M. M., Moklyachuk, M. P., 2012. Interpolation of functionals of stochactic sequanses with stationary increments. Theory Probability and Mathematical Statistics 87, 94-108.

\bibitem{Luz2012b}
Luz, M. M., Moklyachuk, M. P., 2012. Interpolation of functionals of stochastic sequences with stationary increments for observations with noise. Prykladna Statystyka. Aktuarna ta Finansova Matematyka 2, 131-148.

\bibitem{Moklyachuk1994}
Moklyachuk, M. P., 1994. Stochastic autoregressive sequences and minimax interpolation. Theory Probability and Mathematical Statistics 48, 95-103.

\bibitem{Moklyachuk2000}
Moklyachuk, M. P., 2000. Robust procedures in time series analysis. Theory Stochastic Processes 6(3-4), 127-147.

\bibitem{Moklyachuk2001}
Moklyachuk, M. P., 2001. Game theory and convex optimization methods in robust estimation problems. Theory Stochastic Processes 7(1-2), 253-264.

\bibitem{Moklyachuk2008}
Moklyachuk, M. P., 2008. Robust estimations of functionals of stochastic processes. Vydavnycho-Poligrafichny\u\i\ Tsentr, Ky{\"\i}vsky\u\i\ Universytet, Ky{\"\i}v, 320.

\bibitem{Masyutka2006}
Moklyachuk, M.P., Masyutka, O.Yu., 2006. Interpolation of multidimensional stationary sequences. Theory Probability and Mathematical Statistics 73, 125-133.

\bibitem{Masyutka2007}
Moklyachuk, M.P., Masyutka, O.Yu., 2007. On the problem of filtration of vector stationary sequences. Theory Probability and Mathematical Statistics 75, 109-119.

\bibitem{Masyutka2008}
Moklyachuk, M.P., Masyutka, O.Yu., 2008. Minimax prediction problem for multidimensional stationary stochastic sequences. Theory Stochastic Processes 14(3-4),89-103.

\bibitem{Masyutka2011}
Moklyachuk, M.P., Masyutka, O.Yu., 2011. Minimax prediction problem for multidimensional stationary stochastic processes. Communications in Statistics. Theory and Methods 40, 3700-3710.

\bibitem{Masyutka2012}
Moklyachuk, M. P., Masyutka, O.Yu., 2012. Minimax-robust estimation technique for stationary stochastic processes. LAP Lambert Academic Publishing, 296.

\bibitem{Pinsker1954}
Pinsker, M. S., Yaglom, A. M., 1954. On linear extrapolaion of random processes with nth stationary incremens. Doklady Akademii Nauk SSSR 94, 385-388.

\bibitem{Pinsker1955}
Pinsker, M. S., 1955. The theory of curves with nth stationary incremens in Hilber spaces. Izvestiya Akademii Nauk SSSR. Ser. Mat. 19(5), 319-344.


\bibitem{Pshenichnyi}
Pshenichnyi, B.N., 1971. Necessary conditions for an extremum. Pure and Applied mathematics. 4. New York: Marcel Dekker, 1971.

\bibitem{Rozanov}
Rozanov, Yu. A., 1990. Stationary stochastic processes. 2nd rev. ed. Nauka, Moskva, 272.

\bibitem{VastolaPoor}
Vastola, K. S., Poor, H. V., 1983. An analysis of the effects of spectral uncertainty on Wiener filtering. Automatica 28, 289-293.

\bibitem{Wiener}
Wiener, N., 1966. Extrapolation, interpolation and smoothing of stationary time series. Whis wngineering applications. The M. I. T. Press, Massachusetts Institute of Technology, Cambridge, Mass., 163

\bibitem{Yaglom1987a}
Yaglom, A. M., 1987. Correlation theory of stationary and related random functions. Vol. 1: Basic results. Springer Series in Statistics, Springer-Verlag, New York etc., 526.

\bibitem{Yaglom1987b}
Yaglom, A. M., 1987. Correlation theory of stationary and related random functions. Vol. 2: Suplementary notes and references. Springer Series in Statistics, Springer-Verlag, New York etc., 258.

\bibitem{Yaglom1955}
Yaglom, A. M., 1955. Correlation theory of stationary and related random processes with stationary nth increments. Mat. Sbornik 37(1), 141-196.

\bibitem{Yaglom1957}
Yaglom, A. M., 1957. Some clases of random fields in n-dimentional space related with random stationary processes. Teor. Veroyatn. Primen. 2, 292-338







\end{thebibliography}
\end{document}